\documentclass[final,5p,times,twocolumn,numbers]{elsarticle}

\usepackage{amssymb}
\usepackage{amsthm}
\usepackage{amsmath}
\usepackage{tabu} 
\usepackage[T1]{fontenc}
\usepackage{hyperref}
\usepackage{csquotes}
\usepackage{caption}
\usepackage{subcaption}
\usepackage{multirow}
\usepackage{bm}
\usepackage{arydshln}

\newcommand{\vect}[1]{\mathbf{\boldsymbol{#1}}}

\newcommand\restr[2]{{\left.\kern-\nulldelimiterspace#1\right|_{#2}}}

\newcommand{\norm}[1]{\left\lVert#1\right\rVert}

\DeclareMathOperator*{\argmin}{arg\,min}

\DeclareMathOperator*{\opt}{opt}

\makeatletter
\newcommand\xleftrightarrow[2][]{%
  \ext@arrow 9999{\longleftrightarrowfill@}{#1}{#2}}
\newcommand\longleftrightarrowfill@{%
  \arrowfill@\leftarrow\relbar\rightarrow}
\makeatother

\newdefinition{definition}{Definition}

\journal{{}}

\begin{document}

\begin{frontmatter}

\title{A novel reliability-based robust design multi-objective optimization formulation applied in chemical engineering}

\author[1,2]{Gustavo Barbosa Libotte\corref{cor1}}
\ead{gustavolibotte@iprj.uerj.br, glibotte@lncc.br}
\cortext[cor1]{Corresponding author}
\author[3]{Fran S{\'e}rgio Lobato}
\ead{fslobato@ufu.br}
\author[1]{Francisco Duarte Moura Neto}
\ead{fmoura@iprj.uerj.br}
\author[4]{Gustavo Mendes Platt}
\ead{gmplatt@furg.br}

\address[1]{Polytechnic Institute, Rio de Janeiro State University, Nova Friburgo, Brazil}

\address[2]{National Laboratory for Scientific Computing, Petrópolis, Brazil}

\address[3]{Chemical Engineering Faculty, Federal University of Uberl{\^{a}}ndia, Uberl{\^{a}}ndia, Brazil}

\address[4]{School of Chemistry and Food, Federal University of Rio Grande, Santo Ant\^{o}nio da Patrulha, Brazil}

\begin{abstract}
Our fundamental objective is to formulate a multi-objective optimization problem capable of handling robust and reliability-based optimizations, to obtain solutions that satisfy prescribed reliability levels and are least sensitive to external noise. In project optimization, these models play a fundamental role, allowing to obtain parameters and attributes capable of enhancing product performance, reducing costs and operating time. We consider two different approaches capable of quantifying uncertainties during the optimization of mathematical models. In the first, robust optimization, the sensitivity of decision variables in relation to deviations caused by external factors is evaluated. Robust solutions tend to reduce deviations due to possible system changes. The second approach, reliability-based optimization, measures the probability of system failure and obtains model parameters that ensures an established level of reliability. The proposed multi-objective approach, reinforced by a robust strategy, and a guaranteed reliability level, proves to be flexible, while incorporating the uncertainties of the problem under consideration and attending the needs of the project. We tested our novel approach on benchmark and chemical engineering problems, usually treated as deterministic problems, and conclude that the proposed methodology provides a systematic way to evaluate uncertainties, in order to achieve more realistic results, considering external factors.
\end{abstract}




\end{frontmatter}


\section{Introduction} \label{sec:introduction}

Increasingly, the production of new technologies and the creation of engineering solutions are influenced by computer simulations. The mathematical formulation of physical problems is essential to the theoretical development of computational models of products and processes, since it allows the analysis of several factors that can influence the final result, before its effective production. These computer simulations provides a prior analysis of possible adversities, reducing the time required for the project. The optimization of projects has become a fundamental tool in an attempt to enhance the performance of a system or product, while keeping costs low.

Most real problems are subject to multiple and conflicting objectives, which can be related to each other, as pointed out by Babu~et~al.~\cite{bib:babu2005}. The modeling of complex problems using only one objective can be an impractical approach, causing eventual errors to be introduced in the mathematical modeling, when certain simplifications are imposed to meet the chosen modeling paradigm. In contrast, multi-objective optimization models have additional degrees of freedom in relation to those with a single objective, in such a way that there is not a single optimal solution, but a set of Pareto optimal solutions~\cite{bib:collette2004}, any of which is a fair choice. This flexibility allows greater rigor in relation to the description of the physical problem.

Moreover, in practice, several values considered in mathematical models are subject to variations and uncertainties. In addition to possible inconsistencies that may exist during computational modeling, the result of an optimization problem may change due to disregarding the real effect of certain parameters, in order to reduce the computational effort, by assigning deterministic values~\cite{bib:jurecka2007}. A consistent way of considering the uncertainties inherent in a given mathematical model is through robust optimization. In general, robustness can be understood as the property where the performance of the process or technology is minimally sensitive to external factors. Therefore, in robust design optimization, the main focus is on minimizing the effects of variations and uncertainties on the product's performance, while its performance is optimized~\cite{bib:taguchi1999}.

Furthermore, models involving stochastic parameters can present large fluctuations in the optimal results, even if such variations are small due to a high sensitivity of the model on these parameters. In addition, the variability related to manufacturing tolerances and geometric imprecision can affect the project in practice~\cite{bib:cursi2015}. In such cases, the uncertainties of a project can be measured by estimating the probability of failure. Unlike robust optimization, reliability-based design optimization is an approach that takes into account the estimation of the probability of failure in an engineering system, in order to guarantee the satisfaction of the probabilistic constraints at the desired levels, and to obtain the best compromise between cost and security~\cite{bib:haldar2000}.

As these last two modeling artifacts focus on the uncertainties of the design variables and other parameters, it is beneficial to apply robustness and reliability techniques together. Although this requires more computational efforts when compared to deterministic models, this issue is surpassed by the further gain of obtaining optimal feasible solutions that work in real operating conditions, and without the need to eliminate the source of variability.

A growing body of literature has studied reliability-based robust design optimization (RBRDO). Various approaches have been proposed to solve these problems and, according to Rathod~et~al.~\cite{bib:rathod2013}, the existing formulations can be grouped into three categories: (\textit{i})~percentile difference-based RBRDO model; (\textit{ii})~moment-based RBRDO model and; (\textit{iii})~hybrid quality loss functions-based RBRDO model. The main aspects of these formulations are outlined below.

In the percentile difference-based RBRDO model~\cite{bib:du2004}, the percentile performance strategy for assessing both the objective robustness and probabilistic constraints is proposed. An additional objective is included in the original optimization problem, referring to the percentile difference between the right and left tails of the performance function distribution, in order to measure its variability. Minimizing the performance variation can be thought of as shrinking the distance between the percentiles. However, a major limitation of this strategy is the dependence on arbitrary choices of the percentile range, which can significantly affect the results, mainly in terms of robustness.

Based on the quality loss functions proposed by Taguchi~et~al.~\cite{bib:taguchi2004}, the moment-based RBRDO model~\cite{bib:youn2005} combines three types of objectives in order to quantify the robustness of the system. The optimization problem is transformed, aiming at minimizing an objective function composed by the sum of the contributions of the quality loss functions for N-type, S-type, and L-type of the performance characteristic (response of the system). Such functions are expressed in terms of the mean and standard deviation of the performance characteristic, and represent the uncertainties in terms of the deviation from a target value. A serious drawback with this approach, however, is that it requires the definition of weights associated with each term of the objective function. For complex problems, the arbitrariness of such choices can make it difficult to explore non-convex regions in the search space.

In turn, hybrid quality loss functions-based RBRDO model~\cite{bib:yadav2010} is formulated by combining the concepts of goal programming~\cite{bib:schniederjans1995} and quality loss function. In general, this approach is similar to the previous one regarding the definition of the objective function by means of the sum of quality loss functions, but adopts deviational variables, the difference between the targeted value and the actual value calculated for a given performance characteristic, which are classified as desirable and undesirable. This formulation also suffers from a number of pitfalls related to the choice of weights, in addition to the need for adequate normalization or scaling of goals.

Although robustness, reliability and multi-objective optimization have already been used in the literature, the incorporation of all of them together in models of physical or chemical processes still deserves investigation, in view of their pervasive applicability in these fields. The main objective of this work is to propose a novel formulation for optimization problems, taking into account robustness and reliability of the solutions, through a multi-objective approach. We present a more systematic formulation in terms of uncertainties, where maximizing the reliability index is one of the objectives of the optimization problem. The proposed formulation is analyzed by solving benchmark and chemical engineering problems---which are usually investigated in the literature as deterministic problems. The benefits and disadvantages of including uncertainties in optimization problems using the proposed formulation are discussed in terms of choosing the most appropriate optimal value for each project (post-processing solutions) and the behavior of solutions in the face of uncertainties.

\begin{figure}
    \centering
    \includegraphics[width=\columnwidth]{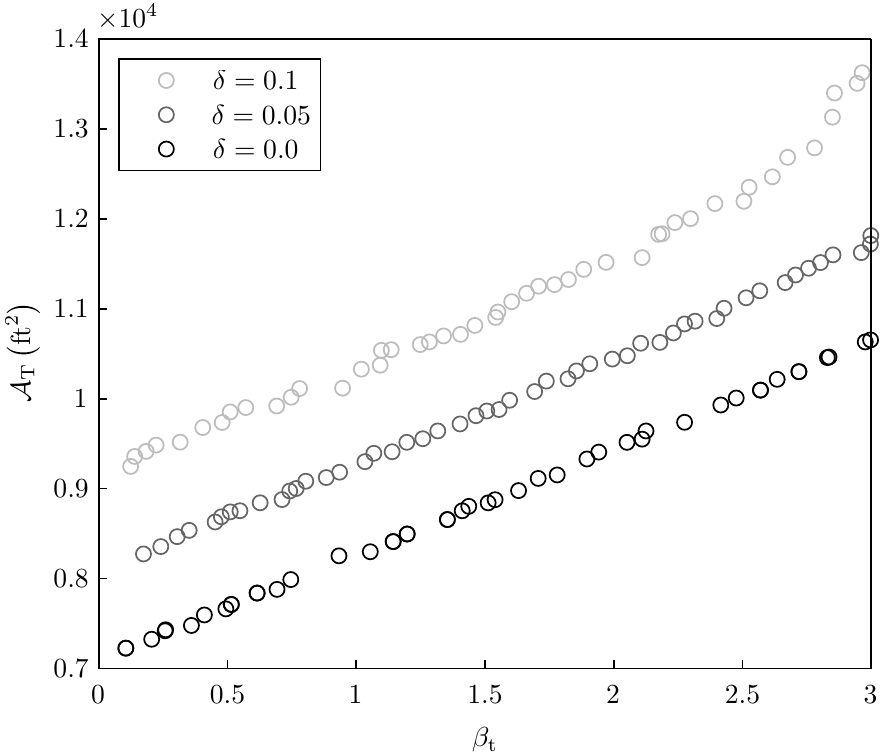}
    \caption{Solutions of the multi-objective heat exchanger network problem with uncertainty, defined by Eq.~(\ref{eq:prob_rede_trocadores_final}), for different levels of robustness.}
    \label{fig:rede_trocadores_pareto}
\end{figure}

Figure~\ref{fig:rede_trocadores_pareto} anticipates some results obtained for the heat exchanger network design problem, which is further discussed in Section~\ref{sec:heat_exchanger_design}. Three sets of non-dominated optimal solutions are shown, obtained from the multi-objective problem, using the metaheuristics Multi-objective Optimization Differential Evolution (MODE)~\cite{bib:lobato2011}. For now, note that the objectives, the total area of the heat exchanger, denoted by $ \mathcal{A}_{\mathrm{T}} $, and the reliability index, given by $ \beta_{\mathrm{t}} $, are conflicting objectives. In this case, more reliable designs of the heat exchanger are related to the increase in area and, therefore, cost. In turn, higher levels of robustness, denoted by $ \delta $, produce Pareto curves that also reflect the increase in the total area, by shifting the curves upwards. Therefore, the novel reliability-based robust design multi-objective formulation provides a tool for the decision maker to strike a balance between robust and reliable solutions, considering the uncertainties associated with the original problem and its model, in order to provide an engineering solution.

The organization of this work is as follows: Section~\ref{sec:multiobjective_thoery} recalls the fundamental concepts of multi-objective optimization and the characteristics of a Pareto set. Section~\ref{sec:robust_optimization} goes on to discuss the specific robustness implementations used. The essential concepts on reliability-based design optimization and our approach to handle non-convex models are shown in Section~\ref{sec:reliability_optimization}. The novel formulation of a multi-objective optimization problem with uncertainty, and a general methodology for evaluating uncertainties, are discussed in Section~\ref{sec:formulation_problem}. Benchmark and chemical engineering multi-objective problems, as well as their reliable and robust optimization, are discussed in Section~\ref{sec:results}. Finally, the conclusions are drawn in Section~\ref{sec:conclusions}.

\section{Optimization and uncertainty}

\subsection{Multi-objective optimization} \label{sec:multiobjective_thoery}

Multi-objective optimization consists in finding a set of points that represents the best balance in relation to optimizing all objectives simultaneously, that is, a collection of solutions that relates the objectives, which are in conflict with each other in most cases~\cite{bib:babu2005}.

Let $ P \subset \mathrm{I\!R}^{n} $ denote a hyperrectangle of all $ \mathbf{x} = \left( x_{1}, \; \dots, \; x_{n} \right)^{\mathrm{T}} \in \mathrm{I\!R}^{n} $, such that $ \mathbf{x}^{\mathrm{inf}} \leq \mathbf{x} \leq \mathbf{x}^{\mathrm{sup}} $, where these inequalities are taken coordinate by coordinate, that is, $ x_{i}^{\mathrm{inf}} \leq x_{i} \leq x_{i}^{\mathrm{sup}} $, with $ x_{i}^{\mathrm{inf}} $ and $ x_{i}^{\mathrm{sup}} $ prescribed, for every $ i = 1, \; \dots, \; n $. In turn, assume that $ \mathbf{f} $ is a function defined in $ P $ with values in $ \mathrm{I\!R}^{m} $, that is, $ \mathbf{f} : P \subset \mathrm{I\!R}^{n} \rightarrow \mathrm{I\!R}^{m} $, with $ m \geq 2 $, where $ \mathbf{f} \left( \mathbf{x} \right) = \left( f_{1} \left( \mathbf{x} \right), \; \dots, \; f_{m} \left( \mathbf{x} \right) \right)^{\mathrm{T}} $ is the \textit{multi-objective function} and $ f_{r} : P \rightarrow \mathrm{I\!R} $, for $ r = 1, \; \dots, \; m $, are the \textit{objectives}. In multi-objective optimization, some of the objectives are to be minimized and others maximized (both are to be optimized). Therefore, let $ m_{1} \geq 0 $ and $ m_{2} \geq 0 $, with $ m = m_{1} + m_{2} $, and let the objectives $ f_{r} $ be ordered such that the first $ m_{1} $ objectives, for $ r = 1, \; \dots, \; m_{1} $, are to be minimized, and the next $ m_{2} $ objectives, for $ r = m_{1} + 1, \; \dots, \; m_{1} + m_{2} $, are to be maximized.

The domain $ P $ of $ \mathbf{f} $ is also called the \textit{design variables on solution space}. In turn, vector $ \mathbf{x} \in P $ is called \textit{decision vector} and its entries are called \textit{decision variables}. In addition, typical problems have additional inequality and equality constraints, which are written in vector form, respectively, as $ \mathbf{g} \left( \mathbf{x} \right) \leq \mathbf{0} $ and $ \mathbf{h} \left( \mathbf{x} \right) = \mathbf{0} $, where $ \mathbf{g} : P \rightarrow \mathrm{I\!R}^{p} $ and $ \mathbf{h} : P \rightarrow \mathrm{I\!R}^{q} $. That is, coordinate-wise, $ g_{i} \left( \mathbf{x} \right) \leq 0 $, for $ i = 1, \; \dots, \; p $, and $ h_{j} \left( \mathbf{x} \right) = 0 $, for $ j = 1, \; \dots, \; q $. The set of solutions that satisfy such constraints, denoted by $ \Omega = \lbrace \mathbf{x} \in \mathrm{I\!R}^{n} \; \vert \; \mathbf{x}^{\mathrm{inf}} \leq \mathbf{x} \leq \mathbf{x}^{\mathrm{sup}}, \; \mathbf{g} \left( \mathbf{x} \right) \leq \mathbf{0}, \; \mathbf{h} \left( \mathbf{x} \right) = \mathbf{0} \rbrace $, is called \textit{feasible set} or \textit{search space}, and its points are termed \textit{feasible solutions}. Another relevant set is the image of the multi-objective function when considering the feasible set, $ \vartheta = \mathrm{Im} \left( \restr{f}{\Omega} \right) \doteq \lbrace \mathbf{y} = \mathbf{f} \left( \mathbf{x} \right) \subset \mathrm{I\!R}^{m}, \; \forall \mathbf{x} \in \Omega \rbrace $, called \textit{objective space}.

Due to the conflict between the objective functions---typically each objective has a different optimizer (either minimizer, if $ r = 1, \; \dots, \; m_{1} $, or maximizer, if $ r = m_{1} + 1, \; \dots, \; m_{1} + m_{2} $) in the feasible set---in general there is no single point capable of optimizing all functions simultaneously. Given a feasible solution, if an algorithm is devised to make a change in it by improving one objective, it may lead to the deterioration of another objective. This raises the question of how the different objectives should be combined to produce a solution in a certain reasonable optimal sense. One reasonable approach is attained by the notion of Pareto optimizer, recalled below after the notion of dominance of feasible solutions~\cite{bib:deb2001}. The concept of dominance establishes a way to compare solutions in the context of multi-objective functions, and most optimization methods search for non-dominated solutions.

A feasible solution $ \mathbf{x}^{1} \in \Omega $ is \textit{dominant} over another solution $ \mathbf{x}^{2} \in \Omega $ (or $ \mathbf{x}^{2} $ is \textit{dominated} by $ \mathbf{x}^{1} $), represented by $ \mathbf{x}^{1} \succ \mathbf{x}^{2} $, if both of the conditions below are met:
\begin{enumerate}
\item Solution $ \mathbf{x}^{1} $ is no worse than $ \mathbf{x}^{2} $ in all objectives, that is, $ f_{r} \left( \mathbf{x}^{1} \right) \leq f_{r} \left( \mathbf{x}^{2} \right) $, for all $ r = 1, \; \dots, \; m_{1} $, and $ f_{r} \left( \mathbf{x}^{1} \right) \geq f_{r} \left( \mathbf{x}^{2} \right) $, for all $ r = m_{1} + 1, \; \dots, \; m $;
\item Solution $ \mathbf{x}^{1} $ is strictly better than $ \mathbf{x}^{2} $ in at least one objective, that is, $ f_{k} \left( \mathbf{x}^{1} \right) < f_{k} \left( \mathbf{x}^{2} \right) $ for some $ k \in \lbrace 1, \; \dots, \; m_{1} \rbrace $ or $ f_{\tilde{k}} \left( \mathbf{x}^{1} \right) > f_{\tilde{k}} \left( \mathbf{x}^{2} \right) $ for some $ \tilde{k} \in \lbrace m_{1} + 1, \; \dots, \; m \rbrace $.
\end{enumerate}

\noindent If any of these conditions is violated, solution $ \mathbf{x}^{1} $ is said not to dominate solution $ \mathbf{x}^{2} $. There are three excluding and exhaustive possibilities: (\textit{i}) $ \mathbf{x}^{1} $ dominates $ \mathbf{x}^{2} $; (\textit{ii}) $ \mathbf{x}^{1} $ is dominated by $ \mathbf{x}^{2} $ and; (\textit{iii}) there is no dominance between $ \mathbf{x}^{1} $ and $ \mathbf{x}^{2} $. The relation $ \succ $ is a strict partial order relation in the set of feasible solutions in $ \Omega $ (irreflexive, asymmetric, and transitive)~\cite{bib:gersting2006}.

A Pareto optimizer is not dominated by any other feasible point in the search space. More formally, a feasible solution $ \mathbf{x}^{*} \in \Omega $ is said to be a \textit{Pareto optimizer} if there is no other feasible vector $ \mathbf{x} \in \Omega $ that dominates $ \mathbf{x}^{*} $, $ \mathbf{x} \succ \mathbf{x}^{*} $, that is, if $ \mathbf{x}^{*} $ is a maximal element of the strict partial order $ \succ $. Denote by $ \mathbf{\Theta} $ the set of Pareto optimizers.

An objective vector $ \mathbf{z}^{*} \in \vartheta $ (that is, exists $ \mathbf{x}^{*} \in \Omega $ at $ \mathbf{z}^{*} = \mathbf{f} \left( \mathbf{x}^{*} \right) $) is said to be a \textit{Pareto optimum} if there is no other objective vector $ \mathbf{z} \in \vartheta $ such that $ z_{i} \leq z_{i}^{*} $ for every $ r = 1, \; \dots, \; m_{1} $, $ z_{i} \geq z_{i}^{*} $ for every $ r = m_{1} + 1, \; \dots, \; m $, and $ z_{k} < z_{k}^{*} $ for some index $ k \in \lbrace 1, \; \dots, \; m_{1} \rbrace $ or $ z_{\tilde{k}} > z_{\tilde{k}}^{*} $ for some index $ \tilde{k} \in \lbrace m_{1} + 1, \; \dots, \; m \rbrace $. Denote by $ \boldsymbol{\mathcal{P}} $ the set of Pareto optima. This is the so-called \textit{Pareto optimal set}, also refereed to as \textit{Pareto front} or \textit{curve} (in the $ m = 2 $ case). Usually, $ \boldsymbol{\mathcal{P}} $ has an infinite number of Pareto optimal values. Note that $ \mathbf{z}^{*} $ is a Pareto optimum if the corresponding feasible vectors, $ \mathbf{x} \in \Omega $ such that $ \mathbf{f} \left( \mathbf{x} \right) = \mathbf{z}^{*} $, are Pareto optimizers. Moreover, $ \mathbf{\Theta} = \mathbf{f}^{-1} \left( \boldsymbol{\mathcal{P}} \right) = \lbrace \mathbf{x} \in \Omega \; \vert \; \mathbf{f} \left( \mathbf{x} \right) \in \boldsymbol{\mathcal{P}} \rbrace $.

The multi-objective optimization problem is denoted by
\begin{align} \label{eq:prob_otim_multi-obj}
\begin{split}
\opt \; &\mathbf{f} \left( \mathbf{x} \right)\\
\textrm{Subject to} \; &\mathbf{g} \left( \mathbf{x} \right) \leq \mathbf{0}\\
\hphantom{\textrm{Subject to}} &\mathbf{h} \left( \mathbf{x} \right) = \mathbf{0}\\
\hphantom{\textrm{Subject to}} &\mathbf{x}^{\mathrm{inf}} \leq \mathbf{x} \leq \mathbf{x}^{\mathrm{sup}} \; ,
\end{split}
\end{align}

\noindent and its solution is $ \boldsymbol{\mathcal{P}} \subset \vartheta $, the Pareto optimal set in the objective space and the corresponding set $ \mathbf{\Theta} $ of Pareto optimizers. If the objectives are all to be minimized, we replace $ \opt $ by $ \min $ in Eq.~(\ref{eq:prob_otim_multi-obj}). Likewise, if all objectives are to be maximized, $ \opt $ is substituted by $ \max $.

Algorithms to solve Eq.~(\ref{eq:prob_otim_multi-obj}) have to deal with the related sets, the Pareto optimal set $ \boldsymbol{\mathcal{P}} $, and the set of Pareto optimizers $ \mathbf{\Theta} $. In the presence of multiple Pareto optimal solutions, that is, when $ \boldsymbol{\mathcal{P}} $ has more than one optimal solution, it does not make sense to choose one result in relation to the others, without any additional information about the problem. In the absence of such information, all of Pareto optimal solutions are equally important. It is thus to the decision maker the task of pinpointing the engineering solution. While it is relevant to find as many Pareto optimal solutions as possible, this translates into two fundamental requirements in the multi-objective optimization algorithms: find a set of values (\textit{i}) as close as possible to the Pareto front (approximation) and; (\textit{ii}) as diverse as possible (covering as much as possible the Pareto front).

\subsection{Robust optimization} \label{sec:robust_optimization}

Sometimes, the objective functions depend sensitively on parameters determined by measurements subject to random noise. This makes for a strong variation on the optimizers. To avoid this, one considers robust solutions, by regularizing the objective function, that is, limiting in some way its overall variability. For this purpose, restrictions on its domain or penalty methods may be adopted. In many practical problems, objective functions are perturbed by random noise, and genetic algorithms~\cite{bib:holland1973,bib:holland1975,bib:dejong1975,bib:goldberg1989} have been widely employed as an effective optimization tool to deal with them. Some of the most recent approaches to tackling robust optimization have emerged in parallel with the progress of techniques related to genetic algorithms.

The ability to solve optimization problems with noise was explored by Fitzpatrick and Grefenstette~\cite{bib:fitzpatrick1988}, who analyzed the effort undertaken by genetic algorithms on such endeavor. Aizawa and Wah~\cite{bib:aizawa1994} developed novel methods to adjust design parameters of genetic algorithms that operate in noisy environments, and B{\"a}ck and Hammel~\cite{bib:back1994} demonstrated that perturbations in the objective function do not significantly influence the performance of evolutionary strategies, as long as the noise level is small enough, compared to the value of the function.

Later, Tsutsui~et~al.~\cite{bib:tsutsui1996} proposed a scheme to search for robust solutions, so that noise, denoted by $ \vect{\delta} \in \mathrm{I\!R}^{n} $, is added to the decision variables $ \mathbf{x} \in \mathrm{I\!R}^{n} $. According to Tsutsui~et~al.~\cite{bib:tsutsui1996}, the evaluation of the effect of noise being added to the output of the system, as previously presented by Hammel and B{\"a}ck~\cite{bib:hammel1994}, is not the most appropriate way to represent the modifications to which the genetic code would be subject to. This approach makes possible to quantify the noise related to each decision variable, making the mathematical model more flexible and capable of describing the physical problem with greater accuracy. Further, Tsutsui and Ghosh~\cite{bib:tsutsui1997} demonstrated the effectiveness of the methodology in obtaining multiple robust solutions.

\subsubsection{Multi-objective robust solution of types I and II} \label{sec:tec_robustez}

Deb and Gupta~\cite{bib:deb2005} extended the formulation presented by Tsutsui and Ghosh~\cite{bib:tsutsui1997}, to compute robust solutions in the multi-objective context, by defining multi-objective robust solution of types I and II. The proposed formulations are essentially derived from the Fundamental Theorem of Genetic Algorithms~\cite{bib:goldberg1989}.

Following Deb and Gupta~\cite{bib:deb2005}, a \textit{multi-objective robust solution of type I} is a feasible Pareto optimizer to the multi-objective optimization problem, Eq.~(\ref{eq:prob_otim_multi-obj}), for $ \mathbf{f}^{\mathrm{eff}} $, subject to $ \mathbf{x} \in \Omega $, such that $ \mathbf{f}^{\mathrm{eff}} \left( \mathbf{x} \right) = \left( f_{1}^{\mathrm{eff}} \left( \mathbf{x} \right), \; \dots, \; f_{m}^{\mathrm{eff}} \left( \mathbf{x} \right) \right)^{\mathrm{T}} $ is the mean effective function, defined with respect to a $ \vect{\delta} $-neighborhood, where
\begin{equation} \label{eq:integral_media_efetiva}
f_{r}^{\mathrm{eff}} \left( \mathbf{x} \right) = \dfrac{1}{\vert \mathcal{B}_{\vect{\delta}} \left( \mathbf{x} \right) \vert} \int\limits_{\vect{\xi} \, \in \, \mathcal{B}_{\vect{\delta}} \left( \mathbf{x} \right)} f_{r} \left( \vect{\xi} \right) \mathrm{d} \vect{\xi} \; ,
\end{equation}

\noindent with $ r = 1, \; \dots, \; m $. $ \mathcal{B}_{\vect{\delta}} \left( \mathbf{x} \right) $ denotes the ball of radius $ \vect{\delta} $ centered at the feasible vector $ \mathbf{x} \in \mathrm{I\!R}^{n} $, and $ \vert \mathcal{B}_{\delta} \left( \mathbf{x} \right) \vert $ is its hypervolume. This definition, also referred to in the literature as effective mean, is appropriate for problems with a closed-form solution to Eq.~(\ref{eq:integral_media_efetiva}). 

For problems in which the search space is more complex, Eq.~(\ref{eq:integral_media_efetiva}) can be approximated by the mean effective value of the objective function, defined by the right-hand side of the equation
\begin{equation*}
f_{i}^{\mathrm{eff}} \left( \mathbf{x} \right) \approx \dfrac{1}{M} \sum\limits_{j \, = \, 1}^{M} f_{i} \left( \vect{\xi}^{j} \right) \; ,
\end{equation*}

\noindent for $ i = 1, \; \dots, \; m $. In this case, $ M $ represents the finite number of random samples taken around $ \mathbf{x} $. Considering the noise vector $ \vect{\delta} $, associated with decision variables $ \mathbf{x} $, points distributed uniformly at random---or respecting some structured scheme---must be arranged in such a way that $ \xi_{i} \in \left[ x_{i} - \delta_{i} x_{i}, \; x_{i} + \delta_{i} x_{i} \right] $, with $ \delta_{i} \in \mathrm{I\!R}_{+} $, for $ i = 1, \; \dots, \; n $.

Based on the previous strategy, Deb and Gupta~\cite{bib:deb2005} also proposed another formulation to quantify the sensitivity of decision variables regarding multi-objective functions, by using a more pragmatic approach, which considers the normalized difference between the values of $ \mathbf{f}^{\mathrm{eff}} $ and $ \mathbf{f} $, in order to provide a threshold capable of quantifying the effect of noise. Hence, a solution is considered robust if such a normalized difference satisfies a prescribed robustness target.

A \textit{multi-objective robust solution of type II} is a feasible Pareto optimizer of the multi-objective optimization problem for $ \mathbf{f} \left( \mathbf{x} \right) $, with feasible set defined by $ \Vert \mathbf{f}^{\mathrm{eff}} \left( \mathbf{x} \right) - \mathbf{f} \left( \mathbf{x} \right) \Vert / \Vert \mathbf{f} \left( \mathbf{x} \right) \Vert \leq \eta $ and $ \mathbf{x} \in \Omega $. Here, the perturbed function $ \mathbf{f}^{\mathrm{eff}} $ is given by Eq.~(\ref{eq:integral_media_efetiva}), or alternatively selecting the worst function value in the objective space among $ M $ samples evaluated within the neighborhood $ \mathcal{B}_{\vect{\delta}} \left( \mathbf{x} \right) $. The parameter $ \eta $ is related to the robustness control, limiting the relative difference between the perturbed and original objective vectors.

\subsubsection{Penalty-based approach} \label{sec:tec_rob_penal}

Mirjalili and Lewis~\cite{bib:mirjalili2016} propose to evaluate the sensitivity of the objective function in the neighborhood of a candidate solution by penalizing the objective function proportionally as it approaches the candidate solution. Unstable points tend to be penalized with greater intensity, according to the fluctuation level of a set of random samples in its neighborhood. The penalty-based robust multi-objective optimization problem is defined as the optimization of $ \mathbf{f}^{\mathrm{P}} \left( \mathbf{x} \right) = \left( f_{1}^{\mathrm{P}} \left( \mathbf{x} \right), \; \dots, \; f_{m}^{\mathrm{P}} \left( \mathbf{x} \right) \right)^{\mathrm{T}}$, subject to $ \mathbf{x} \in \Omega $, with $ f_{r}^{\mathrm{P}} \left( \mathbf{x} \right) = f_{r} \left( \mathbf{x} \right) + P_{r} \left( \mathbf{x} \right) $, for $ 1 \leq r \leq m_{1} $, and $ f_{r}^{\mathrm{P}} \left( \mathbf{x} \right) = f_{r} \left( \mathbf{x} \right) - P_{r} \left( \mathbf{x} \right) $ otherwise. The penalization term is given by
\begin{equation} \label{eq:funcao_penalizacao}
P_{r} \left( \mathbf{x} \right) = \dfrac{\dfrac{1}{M} \sum\limits_{j \, = \, 1}^{M} \left\vert f_{r} \left( \vect{\xi}^{j} \right) - f_{r} \left( \mathbf{x} \right) \right\vert}{\left\vert f_{r} \left( \mathbf{x} \right) \right\vert} \; ,
\end{equation}

\noindent for $ r = 1, \; \dots, \; m $, where $ \vect{\xi} $ denotes random samples around $ \mathbf{x} $, and $ M $ represents the number of random samples.

When small changes in the decision vector cause significant changes in the objective vector, penalty tends to become more pronounced. Consequently, when the objective function is severely penalized for a given candidate solution, it can be considered that this point is very sensitive to external noise and, therefore, not robust, that is, the objective function does not undergo sudden changes when the candidate solution is robust.

\subsection{Reliability-based design optimization} \label{sec:reliability_optimization}

Commonly, repeated observations of natural phenomena produce multiple responses, which can be characterized by the variability in their frequency of occurrence. Such pattern can be pinpointed as an inherent uncertainty of the phenomenon, which is often modeled stochastically. Extending these concepts to engineering systems, where one of the main tasks is to design a system ensuring satisfactory performance, the parameters of an engineering project may present some degree of uncertainty, and therefore must be modeled considering random variables. Essentially, a favorable result is determined in terms of the probability of success, the so-called \textit{reliability} related to the satisfaction of some performance criterion.

For a system which depends on a set of $ n $ continuous normal random variables, given by $ \mathbf{X} = \left( X_{1}, \; \dots, \; X_{n} \right) $, assume that its performance is measured by a scalar $ Y $, the \textit{performance function}, represented by $ Y = g \left( \mathbf{X} \right) $. In turn, the set in $ \mathbf{X} $-space with $ g \left( \mathbf{X} \right) = 0 $ is known as the \textit{limit state function}. When $ Y $ subceed zero, the performance of the system is impaired, that is, it represents states in which the system can no longer fulfill the assignment for which it was designed for. Therefore, the probability of failure is given by
\begin{equation} \label{eq:integral_analise_confiabilidade}
p_{\mathrm{f}} = \mathrm{P} \left[ g \left( \mathbf{X} \right) \leq 0 \right] = \int\limits_{g \left( \mathbf{x} \right) \, \leq \, 0} f_{\mathbf{X}} \left( \mathbf{x} \right) \mathrm{d} \mathbf{x} \; ,
\end{equation}

\noindent where $ f_{\mathbf{X}} \left( \mathbf{x} \right) $ represents the joint probability density function of the random vector $ \mathbf{X} $ and the integration is performed on the failure region. Here, random variables are written in upper case roman letters and their observations in lower case letters.

In order to simplify the computations, it is convenient to perform a linear mapping of the random variables with normal distribution in a set of normalized and independent variables, by transforming the random variable $ \mathbf {X} $ from its original space into a standard normal space $ \mathbf {U} $. This is accomplished by the so-called \textit{Rosenblatt transformation}~\cite{bib:rosenblatt1952}. A particularly simple case arises when the random variables are mutually independent, with $ \mathbf{X} \sim \mathcal{N} \left( \vect{\mu}, \; \vect{\sigma}^{2} \right) $, by means of $ U_{i} = \left( X_{i} - \mu_{i} \right) / \sigma_{i} $, for $ i = 1, \; \dots, \; n $, where the vectors $ \vect{\mu} $ and $ \vect{\sigma} $ are the mean and the standard deviation of the random variables, respectively. The original limit state function $ g \left( \mathbf{X} \right) = 0 $ is mapped into the space $ \mathbf{U} $, where it is represented by $ G \left( \mathbf{U} \right) = 0 $. Thus, the probability of failure, shown in Eq.~(\ref{eq:integral_analise_confiabilidade}), is calculated by the integral
\begin{equation} \label{eq:integral_analise_confiabilidade_transformada}
p_{\mathrm{f}} = \mathrm{P} \left[ G \left( \mathbf{U} \right) \leq 0 \right] = \int\limits_{G \left( \mathbf{u} \right) \, \leq \, 0} \phi \left( \mathbf{u} \right) \mathrm{d} \mathbf{u} \; ,
\end{equation}

\noindent where $ \phi \left( \mathbf{u} \right) $ denotes the standard normal joint probability density function.

Due to the complexity regarding the calculation of a closed-form solution for the integral of Eq.~(\ref{eq:integral_analise_confiabilidade_transformada}), Tu~et~al.~\cite{bib:tu1999} introduced the Performance Measure Approach (PMA), or inverse reliability analysis. In this context, the probability of failure is represented by $ p_{\mathrm{f}} = \Phi \left( -\beta_{\mathrm{t}} \right) $, where $ \Phi $ is the standard normal cumulative distribution function, $ \beta_{\mathrm{t}} \in \mathrm{I\!R}_{+} $ is the \textit{reliability index}, and $ \Phi \left( \beta_{\mathrm{t}} \right) $ denotes the reliability. Note that the larger the reliability index $ \beta_{\mathrm{t}} $ is, the smaller is the probability of failure. PMA determines whether a given design satisfies the probabilistic constraint for a specified target probability of failure, defined in terms of the reliability index, instead of calculating the probability of failure directly. Such formulation consists in minimizing the performance function $ G $, restricted to a hypersphere in $ \mathbf {U} $-space, that is,
\begin{equation} \label{eq:formulacao_PMA}
\begin{split}
\min\limits_{\mathbf{u}} \; &G \left( \mathbf{u} \right)\\
\textrm{Subject to} \; &\norm{\mathbf{u}} = \beta_{\mathrm{t}} \; .
\end{split}
\end{equation}

\noindent Therefore, the value of $ G \left( \mathbf{u} \right) $ is minimized, while the distance between $ \mathbf{u} $ and the origin of the standard normal space remains constant, and equal to the radius $ \beta_{\mathrm{t}} $. In this approach, the hypersphere of radius $ \beta_{\mathrm{t}} $, centered on the origin of the standard normal space $ \mathbf {U} $, is tangent to the limit state function $ G \left( \mathbf{u} \right) = 0 $, and the tangent point is the most probable point of failure, $ \mathbf{u}^{*} $, the solution of Eq.~(\ref{eq:formulacao_PMA}).

Several methods have been proposed to calculate the most probable point of failure using Eq.~(\ref{eq:formulacao_PMA}), among which one mentions Hybrid Modified Chaos Control Method~\cite{bib:meng2015}, Adaptive Modified Chaos Control Method~\cite{bib:li2015}, Relaxed Mean Value Method~\cite{bib:keshtegar2016}, Step Length Adjustment~\cite{bib:yi2016}, Self-adaptive Modified Chaos Control Method~\cite{bib:keshtegar2017a}, Self-adaptive Conjugate Gradient Method~\cite{bib:keshtegar2018b}, Limited Descent-based Mean Value Method~\cite{bib:yaseen2018}, and Adaptive-conjugate Method~\cite{bib:meng2019}. The following section briefly presents the Adaptive Second Order Step Length algorithm~\cite{bib:libotte2020}, employed to solve the probabilistic constraints in this work.

\subsubsection{Adaptive second order step length algorithm} \label{sec:asosl}

The Adaptive Second Order Step Length (ASOSL) algorithm, proposed by Libotte~et~al.~\cite{bib:libotte2020}, is based on the PMA approach, defined by Eq.~(\ref{eq:formulacao_PMA}). The algorithm is a new version of the steepest descent method, but designed for inverse reliability analysis. It uses a second order adaptive step length calculation, without the need to explicitly calculate the Hessian matrix, to accelerate the convergence of the line search. This method is fast and avoids oscillatory behavior.

Let $ G : \mathrm{I\!R}^{n} \rightarrow \mathrm{I\!R} $ be a twice continuously differentiable performance function after transforming the random variables into the standard normal space. Consider an initial guess $ \mathbf{u}^{\left( 0 \right)} \in \mathrm{I\!R}^{n} $ of the most probable point of failure, compute $ G \left( \mathbf{u}^{\left( 0 \right)} \right) $ and the gradient vector $ \mathbf{d}^{\left( 0 \right)} = \left. \nabla G \right\rvert_{\mathbf{u}^{\left( 0 \right)}} $. Then, the computation of the next estimate $ \mathbf{u}^{\left( 1 \right)} $ depends on the step length $ \tau^{\left( 0 \right)} $, given by
\begin{equation} \label{eq:busca_linear_passo_otimo}
\tau^{\left( k \right)} = \argmin\limits_{0 \, < \, t \, \leq \, \bar{t}^{\, \left( k \right)}} \; G \left( \mathbf{u}^{\left( k \right)} - t \mathbf{d}^{\left( k \right)} \right) \; ,
\end{equation}

\noindent with $ \bar{t}^{\, \left( 0 \right)} = 1 $, which is obtained by the backtracking line search procedure~\cite{bib:nocedal2006}.

Next, compute $ \mathbf{u}_{\tau}^{\left( 1 \right)} = \mathbf{u}^{\left( 0 \right)} - \tau^{\left( 0 \right)} \mathbf{d}^{\left( 0 \right)} $. The new approximation of the most probable point of failure is calculated by normalizing $ \mathbf{u}_{\tau}^{\left( 1 \right)} $, $ \mathbf{u}^{\left( 1 \right)} = \beta_{\mathrm{t}} \left( \mathbf{u}_{\tau}^{\left( 1 \right)} / \norm{\mathbf{u}_{\tau}^{\left( 1 \right)}} \right) $, so it belongs to the hypersphere of radius $ \beta_{\mathrm{t}} $. Immediately, one calculates $ G \left( \mathbf{u}^{\left( 1 \right)} \right) $ and $ \mathbf{d}^{\left( 1 \right)} = \left. \nabla G \right\rvert_{\mathbf{u}^{\left( 1 \right)}} $. Next, one calculates the upper bound for the step length line search interval. At first, choose $ \delta_{\eta} $, the search interval extension parameter, and proceed with the iterative process described in what follows. For $ k \geq 1 $, compute $ \bar{t}^{\, \left( k \right)} $ using
\begin{equation*}
\bar{t}^{\, \left( k \right)} = \dfrac{\left( \tau^{\left( k - 1 \right)} \right)^{2} \left( \mathbf{d}^{\left( k - 1 \right)} \right)^{\mathrm{T}} \mathbf{d}^{\left( k - 1 \right)}}{2 \left( G \left( \mathbf{u}^{\left( k \right)} \right) - G \left( \mathbf{u}^{\left( k - 1 \right)} \right) + \tau^{\left( k - 1 \right)} \left( \mathbf{d}^{\left( k - 1 \right)} \right)^{\mathrm{T}} \mathbf{d}^{\left( k - 1 \right)} \right)} \; .
\end{equation*}

\noindent If $ \bar{t}^{\, \left( k \right)} $ is negative (when $ G $ is non-convex), compute a new estimate for $ \bar{t}^{\, \left( k \right)} $, using
\begin{equation*}
\bar{t}^{\, \left( k \right)} = \dfrac{\left( \tau^{\left( k - 1 \right)} + \eta^{\left( k - 1 \right)} \right)^{2} \left( \mathbf{d}^{\left( k - 1 \right)} \right)^{\mathrm{T}} \mathbf{d}^{\left( k - 1 \right)}}{2 \left( G \left( \mathbf{u}^{\left( k \right)} \right) - G \left( \mathbf{u}^{\left( k - 1 \right)} \right) + \left( \tau^{\left( k - 1 \right)} + \eta^{\left( k - 1 \right)} \right) \left( \mathbf{d}^{\left( k - 1 \right)} \right)^{\mathrm{T}} \mathbf{d}^{\left( k - 1 \right)} \right)} \; ,
\end{equation*}

\noindent where $ \eta^{\left( k - 1 \right)} $ is obtained from
\begin{equation*}
\eta^{\left( k - 1 \right)} = \dfrac{G \left( \mathbf{u}^{\left( k - 1 \right)} \right) - G \left( \mathbf{u}^{\left( k \right)} \right) - \tau^{\left( k - 1 \right)} \left( \mathbf{d}^{\left( k - 1 \right)} \right)^{\mathrm{T}} \mathbf{d}^{\left( k - 1 \right)}}{\left( \mathbf{d}^{\left( k - 1 \right)} \right)^{\mathrm{T}} \mathbf{d}^{\left( k - 1 \right)}} + \delta_{\eta} \; .
\end{equation*}

\noindent In this way, it is guaranteed that $ \bar{t}^{\, \left( k \right)} > 0 $ and this value can be used to calculate the optimum step length, $ \tau^{\left( k \right)} $, obtained by solving the line search problem given by Eq.~(\ref{eq:busca_linear_passo_otimo}).

One proceeds to calculate $ \mathbf{u}_{\tau}^{\left( k + 1 \right)} = \mathbf{u}^{\left( k \right)} - \tau^{\left( k \right)} \mathbf{d}^{\left( k \right)} $, and then the new estimate of the most probable point of failure is obtained by normalizing $ \mathbf{u}_{\tau}^{\left( k + 1 \right)} $,
\begin{equation}
\mathbf{u}^{\left( k + 1 \right)} = \beta_{\mathrm{t}} \dfrac{\mathbf{u}_{\tau}^{\left( k + 1 \right)}}{\norm{\mathbf{u}_{\tau}^{\left( k + 1 \right)}}}
\end{equation}

\noindent for $ k \geq 1 $. Next, evaluate $ G \left( \mathbf{u}^{\left( k + 1 \right)} \right) $ and calculate $ \mathbf{d}^{\left( k + 1 \right)} $. Given a pre-established threshold $ \varepsilon > 0 $, sufficiently small, the procedure is halted when $ \mathcal{E}^{\left( k + 1 \right)} < \varepsilon $, with $ \mathcal{E}^{\left( k + 1 \right)} = \norm{\mathbf{u}^{\left( k + 1 \right)} - \mathbf{u}^{\left( k \right)}} $. For further details, see Libotte~et~al.~\cite{bib:libotte2020}.

\section{Formulation of the novel reliability-based robust multi-objective optimization problem} \label{sec:formulation_problem}

Consider the multi-objective optimization problem, shown in Eq.~(\ref{eq:prob_otim_multi-obj}), as the basis for defining the novel reliability-based robust design multi-objective optimization problem. The proposal is that the new formulation, initially deterministic, that is, when stochastic aspects of the model are not taken into account, incorporates such characteristics by taking advantage of the concepts regarding robustness and reliability.

\subsection{Reliability-based robust multi-objective optimization problem}

Assume that, in the multi-objective mathematical model, the design variables can be described by a mix of deterministic and stochastic independent variables. The first type, the deterministic ones, are denoted by $ \mathbf{d} \in \mathrm{I\!R}^{n_{\mathrm{d}}} $ and bounded by $ \mathbf{d}^{\mathrm{inf}} \leq \mathbf{d} \leq \mathbf{d}^{\mathrm{sup}} $. In turn, the stochastic variables are represented by $ \mathbf{X} \in \mathrm{I\!R}^{n_{\mathrm{s}}} $ and defined in the range $ \mathbf{X}^{\mathrm{inf}} \leq \mathbf{X} \leq \mathbf{X}^{\mathrm{sup}} $, with $ n = n_{\mathrm{d}} + n_{\mathrm{s}} $. Essentially, the proposed formulation differs from the deterministic one regarding the evaluation of the uncertainties related to the design variables.

As presented in Section~\ref{sec:robust_optimization}, robust optimization deals with obtaining optimal solutions considering the influence of the sensitivity of decision variables on the value of the objective function. The sensitivity of decision variables must be evaluated both in the objective function and in the constraints of the problem. In turn, considering the uncertainties of a model in terms of its reliability corresponds to transforming the inequality constraints of the problem into probabilistic constraints and, consequently, obtaining the set of feasible solutions to the multi-objective optimization problem in these conditions means that such solutions satisfy a prescribed level of reliability, characterized by a probability of failure, as discussed in Section~\ref{sec:reliability_optimization}.

We recall that the probability of failure is determined by $ \Phi \left( -\beta_{\mathrm{t}} \right) $, where $ \Phi $ is the standard normal cumulative distribution function, and $ \beta_{\mathrm{t}} $ represents the target reliability index. Also, $ \vect{\delta} $ is the vector of noise associated with each coordinate of $ \mathbf{d} $. Thus, let $ \mathbf{F} \left( \mathbf{d}, \; \mathbf{X}, \; \vect{\delta} \right) = \left( F_{1} \left( \mathbf{d}, \; \mathbf{X}, \; \vect{\delta} \right), \; \dots, \; F_{m} \left( \mathbf{d}, \; \mathbf{X}, \; \vect{\delta} \right) \right)^{\mathrm{T}} $. The definition of the novel reliability-based robust design multi-objective optimization problem is given by
\begin{equation} \label{eq:formulacao_multiobjetivo_incertezas}
\begin{split}
\opt \; &\left( \mathbf{F} \left( \mathbf{d}, \; \mathbf{X}, \; \vect{\delta} \right), \; \beta_{\mathrm{t}} \right)\\
\textrm{Subject to} \; &\mathrm{P} \left[ g_{i} \left( \mathbf{d}, \; \mathbf{X}, \; \vect{\delta} \right) \leq 0 \right] \leq \Phi \left( -\beta_{\mathrm{t}} \right)\\
\phantom{Subject to} \; &h_{j} \left( \mathbf{d}, \; \mathbf{X} \right) = 0\\
\phantom{Subject to} \; &\mathbf{d}^{\mathrm{inf}} \leq \mathbf{d} \leq \mathbf{d}^{\mathrm{sup}}\\
\phantom{Subject to} \; &\beta_{\mathrm{t}}^{\mathrm{inf}} \leq \beta_{\mathrm{t}} \leq \beta_{\mathrm{t}}^{\mathrm{sup}} \; ,
\end{split}
\end{equation}

\noindent where, as before, $ i = 1, \; \dots, \; p $, $ j = 1, \; \dots, \; q $, and $ p $ and $ q $ represent, respectively, the number of inequality and equality constraints.

\subsection{General methodology for incorporating robustness and reliability in multi-objective optimization algorithms}

An algorithm to determine the Pareto optimal solutions of a multi-objective optimizaion problem based on metaheuristics is presented next. First, consider an initial population organized in a $ v \times \left( n + 1 \right) $ matrix $ \boldsymbol{\mathcal{P}}^{\left( 0 \right)} = \left( \mathbf{P}_{1}^{\left( 0 \right)}, \; \dots, \; \mathbf{P}_{v}^{\left( 0 \right)} \right)^{\mathrm{T}} $, where each line of $ \boldsymbol{\mathcal{P}}^{\left( 0 \right)} $ represents an individual of the population. In turn, an individual is represented by $ \mathbf{P}_{\ell}^{\left( 0 \right)} = \left( d_{\ell 1}^{\left( 0 \right)}, \; \dots, \; d_{\ell n}^{\left( 0 \right)}, \; \beta_{\ell}^{\left( 0 \right)} \right) $, for $ \ell = 1, \; \dots, \; v $, where the last entry is a reliability index.

In the evaluation of the objective function, the robustness analysis of the candidate solutions must be considered first. For this purpose, the methodology described here considers the effective mean technique, shown in Section~\ref{sec:tec_robustez}. At first, one must define a noise vector $ \vect{\delta} = \left( \delta_{1}, \; \dots, \; \delta_{n}, \; 0 \right) $, where $ \delta_{s} \in \mathrm{I\!R}_{+} $ is the noise level, proportional to the robustness of each decision variable, for $ s = 1, \; \dots, \; n $. The last coordinate of $ \vect{\delta} $ is zero because the reliability index is not subject to noise.

For each fixed $ \ell $, denote by $ \vect{\mathcal{R}}^{\left( k \right)} $ a $ M \times \left( n + 1 \right) $ matrix with $ M $ randomly generated points (following a prescribed probability distribution) in the vicinity of the candidate solution $ \mathbf{P}_{\ell}^{\left( k \right)} $. One such choice would have $ \mathcal{R}_{js}^{\left( k \right)} $ uniformly distributed in an interval of size $ 2 \delta_{s} $ such that
\begin{equation*}
P_{{\ell}s}^{\left( k \right)} \left( 1 - \delta_{s} \right) \leq \mathcal{R}_{js}^{\left( k \right)} \leq P_{{\ell}s}^{\left( k \right)} \left( 1 + \delta_{s} \right) \; ,
\end{equation*}

\noindent for $ s = 1, \; \dots, \; n $ and $ j = 1, \; \dots, \; M $. Note that the last column of $ \vect{\mathcal{R}}^{\left( k \right)} $ is composed of zeros for all $ k $. Therefore, a new function $ \mathbf{F}^{\mathrm{eff}} = \left( F_{1}^{\mathrm{eff}}, \; \dots, \; F_{m}^{\mathrm{eff}} \right)^{\mathrm{T}} $ constructed to assure robustness of the decision variables, is defined as
\begin{equation} \label{eq:funcao_modificada_robustez}
F_{r}^{\mathrm{eff}} \left( \mathbf{d}, \; \mathbf{X}, \; \vect{\delta} \right) = \dfrac{1}{M} \sum\limits_{j \, = \, 1}^{M} F_{r} \left( \mathbf{\mathcal{R}}_{j}, \; \mathbf{X}, \; \vect{\delta} \right) \; ,
\end{equation}

\noindent where $ r = 1, \; \dots, \; m $ and $ \boldsymbol{\mathcal{R}}_{j} $ represents the line $ j $ of $ \vect{\mathcal{R}}^{\left( k \right)} $.

For each point evaluated in Eq.~(\ref{eq:funcao_modificada_robustez}), it is required to check if the probabilistic constraints are satisfied. This verification is performed by executing an inverse reliability analysis method for each probabilistic constraint. In this case, consider an initial estimate $ \mathbf{x}^{\left( 0 \right)} $ for the random variables of the problem, before transforming to standard normal space $ \mathbf{U} $. Let the result of a single execution of the inverse reliability technique be expressed as $ g_{i}^{*} = g_{i} \left( \mathbf{x}^{*} \right) $, for $ i = 1, \; \dots, \; p $, where $ \mathbf{x}^{*} $ represents the most probable point of failure in the $ \mathbf{X} $-space. When a given probabilistic constraint is met, $ g_{i}^{*} > 0 $. Then, a new objective function $ \mathbf{\mathcal{F}} $, which considers the probabilistic constraints penalizing the robust objective function $ \mathbf{F}^{\mathrm{eff}} $, is defined as
\begin{equation} \label{eq:mod_geral_rob_conf}
\mathcal{F}_{r} \left( \mathbf{d}, \; \mathbf{X}, \; \vect{\delta} \right) = F_{r}^{\mathrm{eff}} \left( \mathbf{d}, \; \mathbf{X}, \; \vect{\delta} \right) + \Psi \sum\limits_{i \, = \, 1}^{p} \max \left( -g_{i}^{*}, \; 0 \right) \; ,
\end{equation}

\begin{figure*}[!ht]
    \centering
    \includegraphics[width=\linewidth]{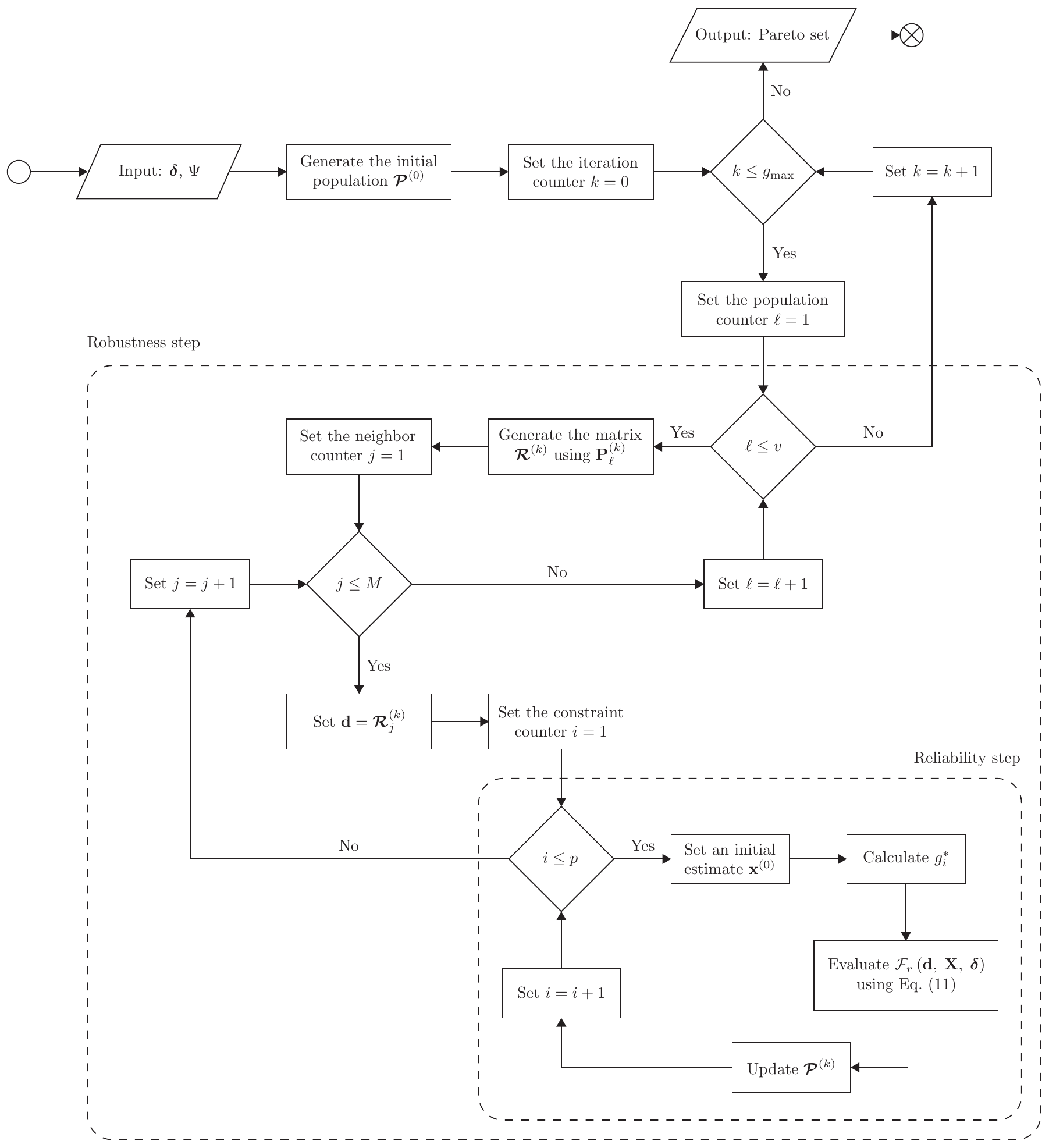}
    \caption{Flowchart of the general methodology for incorporating robustness and reliability in multi-objective optimization problems.}
    \label{fig:flowchart}
\end{figure*}

\noindent where $ \Psi $ is the penalty coefficient, generally adopted as a value large enough for the cases in which $ \mathbf{F} $ is minimized $ \left( r \leq m_{1} \right) $, and $ F_{r}^{\mathrm{eff}} $ is defined by Eq.~(\ref{eq:funcao_modificada_robustez}). When dealing with a maximization problem $ \left( r > m_{1} \right) $, the penalty coefficient is $ -\Psi $. When all probabilistic constraints are satisfied, the contribution related to the reliability analysis is equal to zero in the right-hand side of Eq.~(\ref{eq:mod_geral_rob_conf}) and, therefore, $ \mathcal{F}_{r} \left( \mathbf{d}, \; \mathbf{X}, \; \vect{\delta} \right) = F_{r}^{\mathrm{eff}} \left( \mathbf{d}, \; \mathbf{X}, \; \vect{\delta} \right) $.

Equation~(\ref{eq:mod_geral_rob_conf}) represents a generalization of the function capable of promoting robustness and reliability of the multi-objective optimization problem defined by Eq.~(\ref{eq:formulacao_multiobjetivo_incertezas}). Note that, in the formulation proposed in Eq.~(\ref{eq:mod_geral_rob_conf}), $ \mathbf{\mathcal{F}} $ is the sum of contributions associated with robustness and reliability terms. It turns out that, at the end of the optimization procedure, the Pareto set obtained in the objective space is composed of individuals that are least sensitive to external perturbations and that satisfies the probabilistic constraints. Figure~\ref{fig:flowchart} shows the flowchart with the fundamental elements of the proposed methodology for incorporating robustness and reliability in multi-objective optimization problems.

\section{Results and discussion} \label{sec:results}

In this section the proposed formulation for obtaining robust and reliable solutions in multi-objective optimization problems is evaluated. Initially, a classic benchmark problem is solved, in order to provide an overview of the behavior of the results obtained. Afterwards, some chemical engineering problems are solved. In all cases, the inverse reliability method used to solve the probabilistic constraints is the ASOSL method (Section~\ref{sec:asosl}).

\subsection{Benchmark problem} \label{sec:benchmark_prob}

Initially, let $ f \left( \mathbf{d} \right) = d_{1} + d_{2} $ and consider the constrained single-objective optimization problem, given by
\begin{equation} \label{eq:prob_benchmark_deterministico}
\begin{split}
\min \; &f \left( \mathbf{d} \right)\\
\textrm{Subject to} \; &g_{1} \left( \mathbf{d} \right) = 1 - \dfrac{d_{1}^{2} d_{2}}{20} \leq 0\\
\phantom{Subject to} \; &g_{2} \left( \mathbf{d} \right) = 1 - \dfrac{\left( d_{1} + d_{2} - 5 \right)^{2}}{30} - \dfrac{\left( d_{1} - d_{2} - 12 \right)^{2}}{120} \leq 0\\
\phantom{Subject to} \; &g_{3} \left( \mathbf{d} \right) = 1 - \dfrac{80}{d_{1}^{2} - 8 d_{2} + 5} \leq 0\\
\phantom{Subject to} \; &0 \leq d_{1}, \; d_{2} \leq 10 \; ,
\end{split}
\end{equation}

\noindent In the next four sections, the problem is solved considering (\textit{i}) the classical (deterministic) formulation; (\textit{ii}) the robust approach; (\textit{iii}) the reliability-based approach and; (\textit{iv}) our novel approach, considering robustness and reliability.

\subsubsection{Classical solution}

The solution can be obtained using a numerical optimization method capable of considering inequality constraints, or through a metaheuristic associated with a technique that takes into account such constraints during the optimization procedure. In the second case, Differential Evolution (DE)~\cite{bib:storn1997} can be used in such a way that the objective function is penalized if any of the constraints are violated. The execution of DE for 100 generations, with amplification factor $ F = 0.5 $, crossover probability $ CR = 0.8 $, and $ NP = 50 $ individuals randomly initialized, computes the minimizer of Eq.~(\ref{eq:prob_benchmark_deterministico}), $ \mathbf{d}^{*} = \left( 3.113885, \; 2.062648 \right) $, where $ f \left( \mathbf{d}^{*} \right) = 5.176532 $.

\subsubsection{Robust approach}

\begin{figure}[!ht]
    \centering
    \includegraphics[width=\columnwidth]{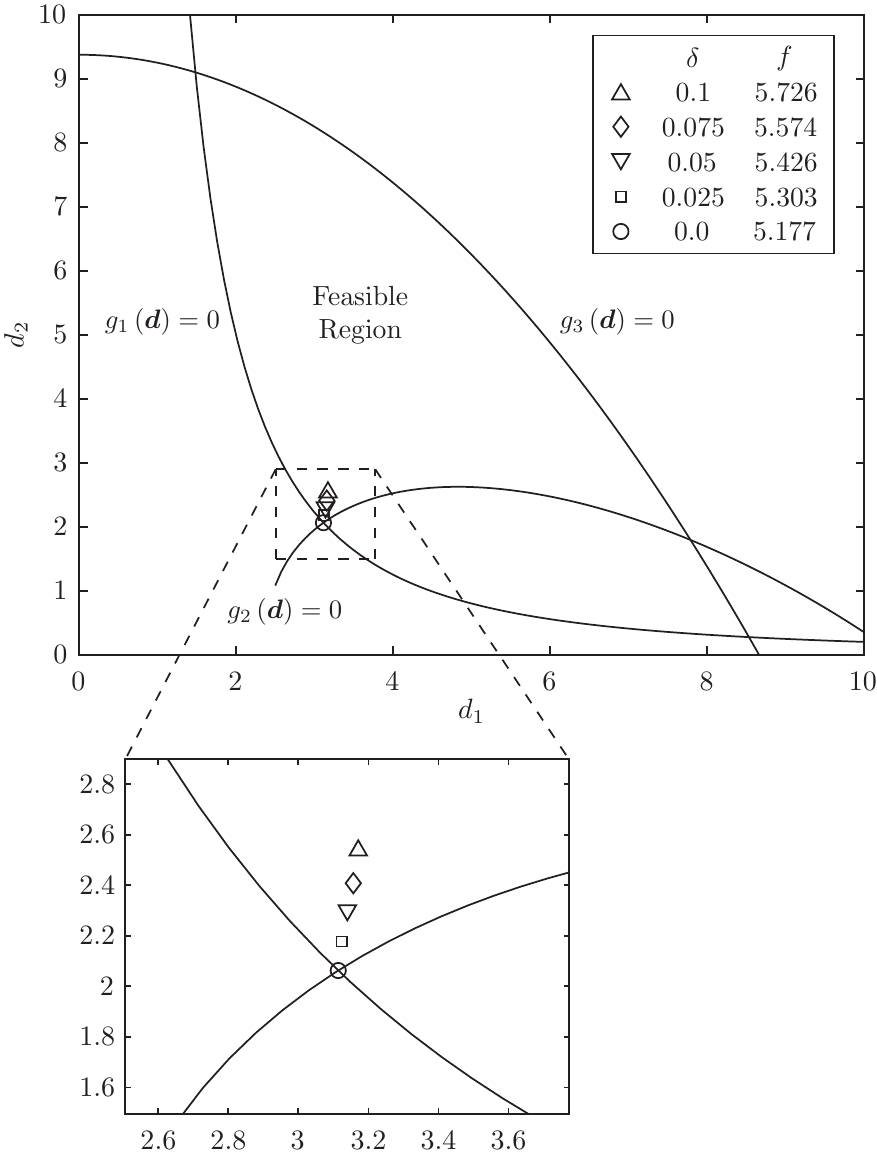}
    \caption{Solution of Eq.~(\ref{eq:prob_benchmark_deterministico}) for different levels of robustness of the design variables. The inserted box shows the value of the objective function corresponding to each level of robustness.}
    \label{fig:robustez_benchmark}
\end{figure}

Now consider the same problem, but in the context of robustness analysis, using the effective mean technique (Section~\ref{sec:tec_robustez}). Again, DE is employed in order to obtain independent solutions for each level $ \delta $ of robustness, using the same parameter values and initialization as in the previous analysis. In the case of the effective mean technique, $ M = 50 $ random points are used, which are generated by the Latin Hypercube strategy~\cite{bib:mencik2016}. The results for values of $ \delta $ ranging from 0\% to 10\% are shown in Fig.~\ref{fig:robustez_benchmark}.

At first, note that the deterministic solution $ \left( \delta = 0 \right) $, represented by the circle in Fig.~\ref{fig:robustez_benchmark}, lies at the intersection point of constraints $ g_{1} = 0 $ and $ g_{2} = 0 $, in which case one says that both are active. Increasing the value of $ \delta $ shows that the optimization method tends to get away from the solution which lies at the intersection point of constraints, that is, solution at which any slight perturbation in the decision vector produce a significant difference in the value of the objective function, making the constraints no longer satisfied. In general, deterministic solutions are sensitive and, therefore, the increase in the noise parameter consists of giving up the optimizer with the best value of the objective function, in favor of obtaining more robust solutions.

\begin{figure}[!hb]
    \centering
    \includegraphics[width=\columnwidth]{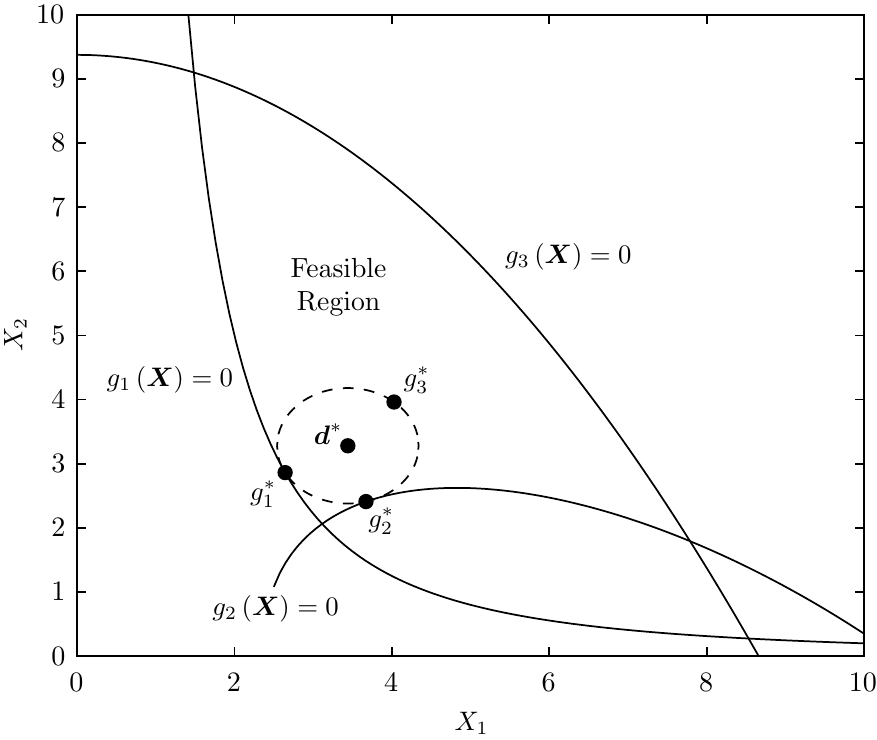}
    \caption{Reliability-based design optimization of the problem given by Eqs.~(\ref{eq:prob_benchmark_conf})--(\ref{eq:restricao_benchmark}) with $ \beta_{\mathrm{t}} = 3 $.}
    \label{fig:confiabilidade_benchmark}
\end{figure}

\subsubsection{A classical reliability approach}

In turn, consider the uncertainties of the problem from the point of view of reliability. In the case of the problem studied, a possible formulation, proposed by Tu~\cite{bib:tu1999b} and adopted in several works in the literature, is given by
\begin{equation} \label{eq:prob_benchmark_conf}
\begin{split}
\min \; &f \left( \mathbf{d} \right)\\
\textrm{Subject to} \; &\mathrm{P} \left[ g_{i} \left( \mathbf{d}, \; \mathbf{X} \right) \leq 0 \right] \leq \Phi \left( -\beta_{\mathrm{t}} \right)\\
\phantom{Subject to} \; &0 \leq d_{1}, \; d_{2} \leq 10 \; ,
\end{split}
\end{equation}

\noindent where $ i = 1, \; 2, \; 3 $ and $ \mathbf{X} $ is a set of independent normal random variables, with $ \vect{\mu} = \left( d_{1}, \; d_{2} \right) $ and $ \vect{\sigma} = \left( 0.3, \; 0.3 \right) $. The target reliability index is $ \beta_{\mathrm{t}} = 3 $, which is equivalent to $ 99.86 \% $ reliability. Constraints are explicitly given by
\begin{subequations} \label{eq:restricao_benchmark}
\begin{equation}
g_{1} \left( \mathbf{X} \right) = 1 - \dfrac{X_{1}^{2} X_{2}}{20}
\end{equation}
\begin{equation}
g_{2} \left( \mathbf{X} \right) = 1 - \dfrac{\left( X_{1} + X_{2} - 5 \right)^{2}}{30} + \dfrac{\left( X_{1} - X_{2} - 12 \right)^{2}}{120}
\end{equation}
\begin{equation}
g_{3} \left( \mathbf{X} \right) = 1 - \dfrac{80}{X_{1}^{2} + 8 X_{2} + 5}\; .
\end{equation}
\end{subequations}

To solve Eqs.~(\ref{eq:prob_benchmark_conf})--(\ref{eq:restricao_benchmark}), DE is employed again, using the same control parameters as in the previous cases (deterministic and robust problems). Probabilistic constraints are solved with ASOSL, where $ \delta_{\eta} = 1 $ and the local search parameters using backtracking are $ \alpha_{\mathrm{b}} = 10^{-4} $ and $ s_{\mathrm{b}} = 0.5 $. Figure~\ref{fig:confiabilidade_benchmark} shows the result obtained in the reliability-based optimization problem. The optimum value reached for the mean of the random variables is $ \mathbf{d}^{*} = \left( 3.440563, \; 3.279963 \right) $, with $ f \left( \mathbf{d}^{*} \right) = 6.720532 $. The dashed circle, shown in Fig.~\ref{fig:confiabilidade_benchmark}, represents the hypersphere of radius $ \beta_{\mathrm{t}} $ (in this case, a circle).

Points $ g_{i}^{*} $, for $ i = 1, \; 2, \; 3 $, represent the values of the probabilistic constraints evaluated at the most probable point of failure, which are obtained during the executions of ASOSL. Note that only constraints $ g_{1} $ and $ g_{2} $ are active, with $ g_{1} \left( 2.6435, \; 2.8619 \right) \approx 0 $ and $ g_{2} \left( 3.6709, \; 2.4099 \right) \approx 0 $. In the case of the third constraint, $ g_{3} \left( 4.0275, \; 3.9623 \right) \approx 0.5118 $ and, therefore, $ g_{3}^{*} > 0 $, which means that all probabilistic constraints are satisfied. In this analysis, it is clear that reliable solutions are always in the center of the circle of radius $ \beta_{\mathrm{t}} $ (for problems in the plane), which tend to move further away from the deterministic minimizer, as the value of $ \beta_{\mathrm{t}} $ is increased. This provides more reliable solutions, but at the expense of increasing the optimal value (which is to be minimized), as is also the case in robust optimization.

\subsubsection{Our novel reliability-based robust approach} \label{sec:benchmark_novel_approach}

Finally, consider the formulation of the reliability-based robust design multi-objective optimization problem, as proposed in Eq.~(\ref{eq:formulacao_multiobjetivo_incertezas}). The problem is formulated as
\begin{equation} \label{eq:multiobjetivo_incertezas_benchmark}
\begin{split}
\opt \; &\left( f \left( \mathbf{d}, \; \vect{\delta} \right), \; \beta_{\mathrm{t}} \right)\\
\textrm{Subject to} \; &\mathrm{P} \left[ g_{i} \left( \mathbf{d}, \; \mathbf{X}, \; \vect{\delta} \right) \leq 0 \right] \leq \Phi \left( -\beta_{\mathrm{t}} \right)\\
\phantom{Subject to} \; &0 \leq d_{1}, \; d_{2} \leq 10\\
\phantom{Subject to} \; &1 \leq \beta_{\mathrm{t}} \leq 3 \; .
\end{split}
\end{equation}

\noindent where $ i = 1, \; 2, \; 3 $ and the expressions for constraints are given by Eq.~(\ref{eq:restricao_benchmark}). Here, $ f $ is to be minimized and $ \beta_{\mathrm{t}} $ is to be maximized. As in the previous case, $ \mathbf{X} $ is a set of independent normal random variables, with $ \vect{\mu} = \left( d_{1}, \; d_{2} \right) $ and $ \vect{\sigma} = \left( 0.3, \; 0.3 \right) $. In this particular case, we choose the reliability index between 1 and 3, which results is a reliability between $ \Phi \left( 1 \right) \approx 0.8413 $ and $ \Phi \left( 3 \right) \approx 0.9987 $.

\begin{figure}[!ht]
    \centering
    \includegraphics[width=\columnwidth]{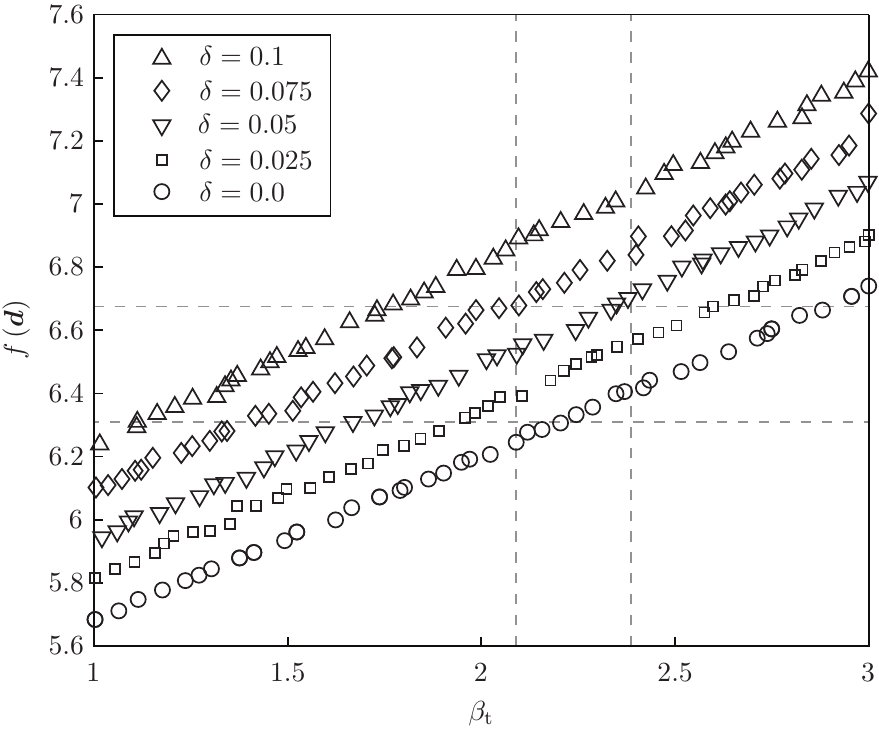}
    \caption{Sets of optimal values obtained for different levels of robustness, in the benchmark problem given by Eq.~(\ref{eq:multiobjetivo_incertezas_benchmark}). The dashed gray lines support comparison of results.}
    \label{fig:robustez_conf_benchmark}
\end{figure}

The MODE~\cite{bib:lobato2011} algorithm is employed to obtain the Pareto set, with amplification factor $ F = 0.5 $, crossover probability $ CR = 0.8 $, reduction parameter $ r = 0.9 $ and number of pseudo-fronts $ R = 10 $. The algorithm runs for $ 500 $ generations for each level of robustness, with $ NP = 50 $ individuals. In the case of probabilistic constraints, ASOSL is used with the same parameters as in the previous case. The computed approximations of the Pareto curves, at different levels of robustness, are shown in Fig.~\ref{fig:robustez_conf_benchmark}.

Analyzing the profiles obtained, the particular contribution of both approaches to uncertainty is clear: maximizing reliability, in terms of $ \beta_{\mathrm{t}} $, represents an increase of the optimal value of $ f $, as the reliability index also increases. In the case of robustness analysis, the profile is shifted, also affecting $ f $, in relation to the Pareto curves obtained when $ \delta $ is smaller. These results are in line with those presented in Figs.~\ref{fig:robustez_benchmark} and \ref{fig:confiabilidade_benchmark}, when the calculated results show a gradual displacement from the deterministic minimizer, but for different reasons. These facts, linked to the economic interests of the project, are fundamental in the choice of solutions to be implemented in practice.

\subsection{Heat exchanger network design} \label{sec:heat_exchanger_design}

\subsubsection{Classical approach and formulation of the problem}

\begin{figure}[!ht]
    \centering
    \includegraphics[width=\columnwidth]{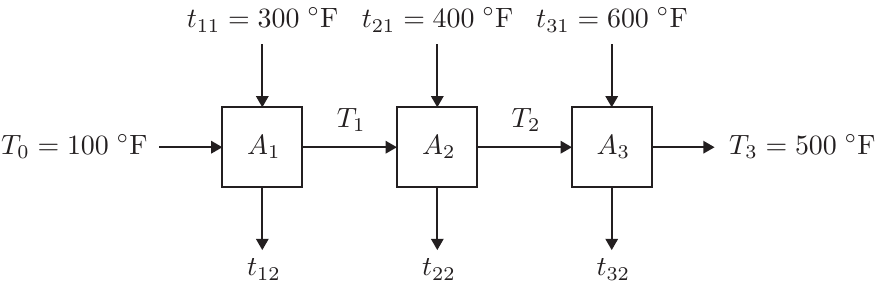}
    \caption{Schematic of the three stage heat exchanger network~\cite{bib:avriel1971}.}
    \label{fig:rede_trocador_calor}
\end{figure}

In the design of a three-stage heat exchanger network in series, as shown in Fig.~\ref{fig:rede_trocador_calor}, constants $ t_{11} $, $ t_{21} $ and $ t_{31} $ represent the temperatures of the hot fluids entering the heat exchangers, and the overall heat transfer coefficients are given by $ U_{1} $, $ U_{2} $ and $ U_{3} $. A fluid, with a given flow rate $ W $ and specific heat $ C_{p} $, is heated from temperature $ T_{0} $ to $ T_{3} $, by passing through three heat exchangers. $ T_{0} $ is the initial temperature of the fluid being heated, $ T_{i} $ is the temperature of the fluid at each heating stage, and $ t_{i2} $, for $ i = 1, \; 2, \; 3 $, are, respectively, the temperatures of the hot fluid leaving the heat exchangers. In each heat exchanger, the cold stream is heated by a hot fluid with the same flow rate $ W $ and specific heat $ C_{p} $~\cite{bib:avriel1971}.

The rate of heat transferred to the cold fluid is less than or equal to the rate of heat lost by the hot stream. Three heat balances express this fact: $ W C_{p} \left( T_{i} - T_{i - 1} \right) \leq W C_{p} \left( t_{i1} - t_{i2} \right) $, for $ i = 1, 2, 3 $. Such expressions can be written as
\begin{subequations} \label{eq:restricoes_heat_exchanger}
\begin{equation} \label{eq:balancos_calor}
T_{i} + t_{i2} \leq t_{i1} + T_{i - 1} \; .
\end{equation}
\noindent Similarly, the heat transfer inequalities for each stage of the network are given by
\begin{equation} \label{eq:inequacoes_calor}
W C_{p} \left( T_{i} - T_{i - 1} \right) \leq U_{i} A_{i} \left( t_{i2} - T_{i - 1} \right) \; ,
\end{equation}
\end{subequations}

\noindent where the overall heat transfer coefficients are given by $ \mathbf{U} = \left( 120, \; 80, \; 40 \right) \; \mathrm{Btu / \left( ft^{2} \, h \, ^{\circ} F \right)} $, and $ W C_{p} = 10^{5} \; \mathrm{Btu / \left( h \, ^{\circ} F \right)} $.

Considering the constraints defined by Eq.~(\ref{eq:restricoes_heat_exchanger}), the objective is to minimize the sum of the heat transfer areas of the three exchangers, given by $ \mathcal{A}_{\mathrm{T}} \left( A_{1}, \; A_{2}, \; A_{3} \right) = A_{1} + A_{2} + A_{3} $. Therefore, the deterministic single-objective optimization problem of the heat exchanger network design is formulated as
\begin{equation*}
\begin{split}
\min \; &\mathcal{A}_{\mathrm{T}} \left( A_{1}, \; A_{2}, \; A_{3} \right)\\
\textrm{Subject to} \; &\dfrac{1}{400} \left( T_{1} + t_{12} \right) - 1 \leq 0\\
\phantom{Subject to} \; &\dfrac{1}{400} \left( T_{2} + t_{22} - T_{1} \right) - 1 \leq 0\\
\phantom{Subject to} \; &\dfrac{1}{100} \left(t_{32} - T_{2} \right) - 1 \leq 0\\
\phantom{Subject to} \; &A_{1} \left( 100 - t_{12} \right) + \dfrac{2500}{3} T_{1} - \dfrac{250000}{3} \leq 0\\
\phantom{Subject to} \; &A_{2} \left( T_{1} - t_{22} \right) - 1250 T_{1} + 1250 T_{2} \leq 0\\
\phantom{Subject to} \; &A_{3} \left( T_{2} - t_{32} \right) - 2500 T_{2} + 1.25 \times 10^{6} \leq 0\\
\phantom{Subject to} \; &10^{2} \leq A_{1} \leq 10^{4}\\
\phantom{Subject to} \; &10^{3} \leq A_{2}, \; A_{3} \leq 10^{4}\\
\phantom{Subject to} \; &10 \leq T_{1}, \; T_{2}, \; t_{12}, \; t_{22}, \; t_{32} \leq 10^{3} \; .
\end{split}
\end{equation*}

\noindent The minimum total area required for the operation of this system is $ \mathcal{A}_{\mathrm{T}} \approx 7049.25 \; \mathrm{ft^{2}} $, with $ A_{1} = 579.31 \; \mathrm{ft^{2}} $, $ A_{2} = 1359.97 \; \mathrm{ft^{2}} $ and $ A_{3} = 5109.97 \; \mathrm{ft^{2}} $. This problem was also studied by Andrei~\cite{bib:andrei2013} and others.

\subsubsection{Incorporating reliability}

Originally, the problem has eight decision variables. However, Angira and Babu~\cite{bib:angira2002} show that inequality constraints, related to Eq.~(\ref{eq:balancos_calor}), can be replaced by equality constraints, without changing the behavior of the system. Rewriting the terms of the temperatures of the hot fluid entering the heat exchangers, $ t_{11} $, $ t_{21} $ and $ t_{31} $, in terms of $ T_{1} $ and $ T_{2} $, and replacing the resultant expressions in the constraints defined by Eq.~(\ref{eq:inequacoes_calor}), the problem now has five decision variables, namely, the area of the heat exchangers, $ A_{1} $, $ A_{2} $ and $ A_{3} $, and temperatures $ T_{1} $ and $ T_{2} $. Thus, the new constraints of the problem are as follows:
\begin{subequations} \label{eq:nova_restricao_rede_trocadores}
\begin{equation} \label{eq:nova_restricao_rede_trocadores_1}
A_{1} T_{1} + \dfrac{2500}{3} T_{1} - 300 A_{1} - \dfrac{250000}{3} \leq 0
\end{equation}
\begin{equation}
T_{2} A_{2} - 400 A_{2} - 1250 \left( T_{1} - T_{2} \right) \leq 0
\end{equation}
\begin{equation} \label{eq:nova_restricao_rede_trocadores_3}
1.25 \times 10^{6} - 2500 T_{2} - 100 A_{3} \leq 0 \; .
\end{equation}
\end{subequations}

For reliability-based optimization, Lobato~et~al.~\cite{bib:lobato2019} suggest that model uncertainties are associated with the constant terms in the inequalities, so that they are treated as independent normal random variables $ \mathbf{X} = \left( X_{1}, \; \dots, \; X_{8} \right) $, with mean equal to the constants in Eq.~(\ref{eq:nova_restricao_rede_trocadores}) and standard deviation as in Table~\ref{tab:prop_estat_rede_trocadores_calor}. In this case, the equivalent inequalities of the deterministic problem, given by Eqs.~(\ref{eq:nova_restricao_rede_trocadores}), can be transformed into probabilistic constraints.

Let $ \mathbf{d} = \left( A_{1}, \; A_{2}, \; A_{3}, \; T_{1}, \; T_{2} \right) $ be the deterministic decision vector and $ \beta_{\mathrm{t}} $ the target reliability index. The heat exchanger network model considering the reliability of solutions can be modeled as
\begin{equation*}
\begin{split}
\opt \; &\left( \mathcal{A}_{\mathrm{T}} \left( d_{1}, \; d_{2}, \; d_{3} \right), \; \beta_{\mathrm{t}} \right)\\
\textrm{Subject to} \; &\mathrm{P} \left[ g_{i} \left( \mathbf{d}, \; \mathbf{X} \right) \leq 0 \right] \leq \Phi \left( -\beta_{\mathrm{t}} \right)\\
\phantom{Subject to} \; &10^{2} \leq d_{1} \leq 10^{4}\\
\phantom{Subject to} \; &10^{3} \leq d_{2}, \; d_{3} \leq 10^{4}\\
\phantom{Subject to} \; &10 \leq d_{4}, \; d_{5} \leq 10^{3}\\
\phantom{Subject to} \; &0.1 \leq \beta_{\mathrm{t}} \leq 3 \; ,
\end{split}
\end{equation*}

\noindent for $ i = 1, \; 2, \; 3 $, where the performance functions are given by the following expressions:
\begin{subequations}
\begin{equation*}
g_{1} \left( \mathbf{d}, \; \mathbf{X} \right) =  d_{1} d_{4} + X_{1} d_{4} - X_{2} d_{1} - X_{3}
\end{equation*}
\begin{equation*}
g_{2} \left( \mathbf{d}, \; \mathbf{X} \right) = d_{5} d_{2} - X_{4} d_{2} - X_{5} \left( d_{4} - d_{5} \right)
\end{equation*}
\begin{equation*}
g_{3} \left( \mathbf{d}, \; \mathbf{X} \right) = X_{6} - X_{7} d_{5} - X_{8} d_{3} \; .
\end{equation*}
\end{subequations}

\begin{table}[!ht]
\centering
\caption{Mean $ \left( \vect{\mu} \right) $ and standard deviation $ \left( \vect{\sigma} \right) $ of the random variables of the heat exchanger network problem. The standard deviation is defined according to the coefficient of variation, which is equal to 0.05.}
\label{tab:prop_estat_rede_trocadores_calor}
{\setlength{\tabulinesep}{1.4mm}
\setlength{\tabcolsep}{3pt}
\begin{tabu}{lcccccccccc}
\hline\hline
& $ X_{1} $ & $ X_{2} $ & $ X_{3} $ & $ X_{4} $ & $ X_{5} $ & $ X_{6} $ & $ X_{7} $ & $ X_{8} $ \\ \hline
$ \vect{\mu} $ & $ \dfrac{2500}{3} $ & 300 & $ \dfrac{250000}{3} $ & 400 & 1250 & $ 1.25 \times 10^{6} $ & 2500 & 100 \\
$ \vect{\sigma} $ & $ \dfrac{125}{3} $ & 15 & $ \dfrac{12500}{3} $ & 20 & $ \dfrac{125}{2} $ & 62500 & 125 & 5 \\ \hline\hline
\end{tabu}}
\end{table}

\subsubsection{Reliability and robustness analysis} \label{sec:heat_exchanger_novel_approach}

In the case of robustness analysis, uncertainties are considered as geometric variations referring to the areas of each heat exchanger ($ A_{1} $, $ A_{2} $ and $ A_{3} $) and differences related to the temperatures measured during fluid flow between the heat exchangers ($ T_{1} $ and $ T_{2} $). Thus, the noise vector $ \vect{\delta} $ determines the permissible range of variation of each of the components of $ \mathbf{d} $. Therefore, the reliability-based robust multi-objective optimization problem of the heat exchanger network design can be formulated as
\begin{equation} \label{eq:prob_rede_trocadores_final}
\begin{split}
\opt \; &\left( \mathcal{A}_{\mathrm{T}} \left( \mathbf{d}, \; \vect{\delta} \right), \; \beta_{\mathrm{t}} \right)\\
\textrm{Subject to} \; &\mathrm{P} \left[ g_{i} \left( \mathbf{d}, \; \mathbf{X}, \; \vect{\delta} \right) \leq 0 \right] \leq \Phi \left( -\beta_{\mathrm{t}} \right)\\
\phantom{Subject to} \; &10^{2} \leq d_{1} \leq 10^{4}\\
\phantom{Subject to} \; &10^{3} \leq d_{2}, \; d_{3} \leq 10^{4}\\
\phantom{Subject to} \; &10 \leq d_{4}, \; d_{5} \leq 10^{3}\\
\phantom{Subject to} \; &0.1 \leq \beta_{\mathrm{t}} \leq 3 \; .
\end{split}
\end{equation}

In order to obtain the solutions to this problem, MODE and ASOSL are adopted, using the same values assigned in the previous problem (Section~\ref{sec:benchmark_prob}) for the control parameters. Figure~\ref{fig:rede_trocadores_pareto} shows the Pareto set obtained considering different levels of robustness, which are calculated by the effective mean strategy, using $ M = 50 $ random points (again with Latin Hypercube sampling). In turn, the reliability ranges from $ \Phi \left( 0.1 \right) \approx 0.5398 $ to $ \Phi \left( 3 \right) \approx 0.9987 $.

\subsubsection{Simulations}

Initially, one may note that the result of the deterministic problem, obtained in the case where $ \delta = 0 $ in Fig.~\ref{fig:rede_trocadores_pareto}, is equivalent to the lowest level of reliability, in the proposed formulation. In general, problems of this type express the total area of the heat exchanger as a function of its performance and operating cost. Therefore, higher costs (as a consequence of the unavoidable increase in the total area of the heat exchanger network) are expected from optimal values whose reliability is maximized, to the detriment of the deterministic result. In turn, when the desired level of robustness for the results is increased, by increasing $ \delta $, the set of dominant solutions is shifted, showing an even more pronounced increase in relation to the total area.

The nearly linear profile obtained for each set of dominant solutions, in the interval in which $ \beta_{\mathrm{t}} $ is evaluated, may be due to the linear behavior of the problem, both in relation to the objective function and the probabilistic constraints. However, the Pareto sets never cross each other. This fact corroborates the consistency of the proposed formulation, as well as the numerical methods used to obtain the results, since Pareto sets with different levels of robustness must not cross each other, making the choice of the optimizer (post-processing step) to be determined considering independent levels of robustness and reliability.

Table~\ref{tab:resultados_extremos_rede_trocadores} shows the values of the decision variables referring to three different optimizers, for each value of $ \vect{\delta} $. The chosen minimizers have different levels of compromise between the objectives: (\textit{i}) solutions in which the total area is minimized; (\textit{ii}) the reliability index is maximized and; (\textit{iii}) solutions with intermediate compromise in relation to both.

\begin{table}[!ht]
\centering
\caption{Several optimizers of the heat exchanger network problem, considering different levels of compromise between the objectives. For each robustness level, the values are displayed in relation to $ \left( A_{1}, \; A_{2}, \; A_{3}, \; T_{1}, \; T_{2} \right)^{\mathrm{T}} $ (above the dashed line), and $ \left( \mathcal{A}_{\mathrm{T}}, \; \beta_{\mathrm{t}} \right)^{\mathrm{T}} $ (below the dashed line).}
\label{tab:resultados_extremos_rede_trocadores}
{\setlength{\tabulinesep}{1.6mm}
\setlength{\tabcolsep}{2.4mm}
\begin{tabu}{cccc}
\hline\hline
\multirow{2}{*}{Robustness} & \multicolumn{3}{c}{Compromise} \\ \cline{2-4}
 & $ \mathcal{A}_{\mathrm{T}} $ & Intermediate & $ \beta_{\mathrm{t}} $ \\ \hline
\multirow{7}{*}{$ \delta = 0.0 $} & $ 725.64 \; \mathrm{ft^{2}} $ & $ 595.07 \; \mathrm{ft^{2}} $ & $ 551.11 \; \mathrm{ft^{2}} $ \\
 & $ 1270.47 \; \mathrm{ft^{2}} $ & $ 1260.02 \; \mathrm{ft^{2}} $ & $ 1279.83 \; \mathrm{ft^{2}} $ \\
 & $ 5229.70 \; \mathrm{ft^{2}} $ & $ 6986.06 \; \mathrm{ft^{2}} $ & $ 8822.10 \; \mathrm{ft^{2}} $ \\
 & $ 190.93 \; \mathrm{^{\circ} F} $ & $ 169.00 \; \mathrm{^{\circ} F} $ & $ 153.48 \; \mathrm{^{\circ} F} $ \\
 & $ 295.11 \; \mathrm{^{\circ} F} $ & $ 269.24 \; \mathrm{^{\circ} F} $ & $ 246.37 \; \mathrm{^{\circ} F} $ \\ \cdashline{2-4}
 & $ 7225.81 \; \mathrm{ft^{2}} $ & $ 8841.15 \; \mathrm{ft^{2}} $ & $ 10653.04 \; \mathrm{ft^{2}} $ \\
 & 0.11 & 1.51 & 3.00 \\ \hline
\multirow{7}{*}{$ \delta = 0.05 $} & $ 467.80 \; \mathrm{ft^{2}} $ & $ 431.84 \; \mathrm{ft^{2}} $ & $ 368.17 \; \mathrm{ft^{2}} $ \\
 & $ 1482.09 \; \mathrm{ft^{2}} $ & $ 1380.90 \; \mathrm{ft^{2}} $ & $ 1497.34 \; \mathrm{ft^{2}} $ \\
 & $ 6323.53 \; \mathrm{ft^{2}} $ & $ 8065.84 \; \mathrm{ft^{2}} $ & $ 9948.55 \; \mathrm{ft^{2}} $ \\
 & $ 155.83 \; \mathrm{^{\circ} F} $ & $ 146.84 \; \mathrm{^{\circ} F} $ & $ 129.49 \; \mathrm{^{\circ} F} $ \\
 & $ 270.75 \; \mathrm{^{\circ} F} $ & $ 247.02 \; \mathrm{^{\circ} F} $ & $ 227.35 \; \mathrm{^{\circ} F} $ \\ \cdashline{2-4}
 & $ 8273.42 \; \mathrm{ft^{2}} $ & $ 9878.58 \; \mathrm{ft^{2}} $ & $ 11814.06 \; \mathrm{ft^{2}} $ \\
 & 0.17 & 1.55 & 3.00 \\ \hline
\multirow{7}{*}{$ \delta = 0.1 $} & $ 141.83 \; \mathrm{ft^{2}} $ & $ 464.59 \; \mathrm{ft^{2}} $ & $ 458.75 \; \mathrm{ft^{2}} $ \\
 & $ 1643.37 \; \mathrm{ft^{2}} $ & $ 1756.62 \; \mathrm{ft^{2}} $ & $ 3178.67 \; \mathrm{ft^{2}} $ \\
 & $ 7461.87 \; \mathrm{ft^{2}} $ & $ 8857.41 \; \mathrm{ft^{2}} $ & $ 9989.32 \; \mathrm{ft^{2}} $ \\
 & $ 113.95 \; \mathrm{^{\circ} F} $ & $ 135.72 \; \mathrm{^{\circ} F} $ & $ 113.86 \; \mathrm{^{\circ} F} $ \\
 & $ 242.45 \; \mathrm{^{\circ} F} $ & $ 244.91 \; \mathrm{^{\circ} F} $ & $ 246.04 \; \mathrm{^{\circ} F} $ \\ \cdashline{2-4}
 & $ 9247.07 \; \mathrm{ft^{2}} $ & $ 11078.62 \; \mathrm{ft^{2}} $ & $ 13626.74 \; \mathrm{ft^{2}} $ \\
 & 0.12 & 1.60 & 2.97 \\ \hline\hline
\end{tabu}}
\end{table}

Note that, for a given level of robustness, it is not possible to verify any obvious linear relationship with respect to the increase in the total area of the heat exchanger network, in detriment to the nearly linear increase observed in Fig.~\ref{fig:rede_trocadores_pareto}. The same goes for temperatures $ T_{1} $ and $ T_{2} $. On the other hand, analyzing each of the minimizers in Table~\ref{tab:resultados_extremos_rede_trocadores}, it is clear that the biggest contribution to the heating of the fluid up to $ 500 \; \mathrm{^{\circ} F} $ is given in the third stage, which corroborates the values of $ A_{3} $ in all reported results, greater than the other heat exchangers in the network.

\subsection{Reactor network design} \label{sec:reactor_network_design}

\begin{figure}[!ht]
    \centering
    \includegraphics[width=\columnwidth]{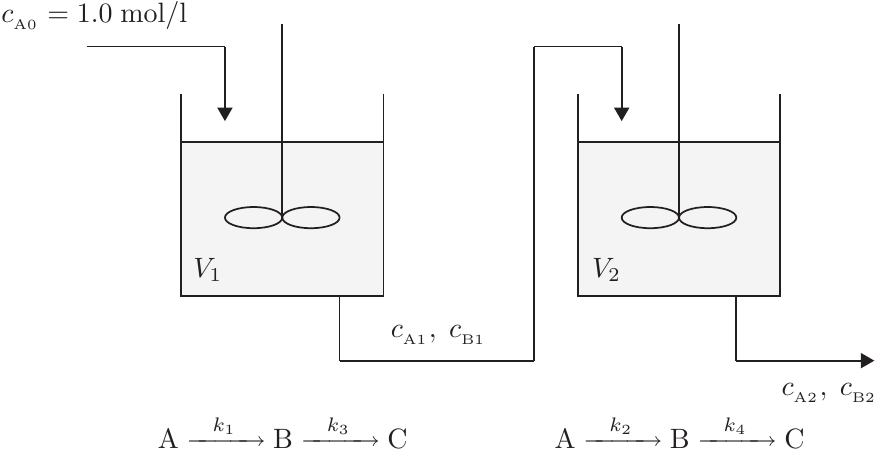}
    \caption{Reactor network scheme~\cite{bib:smith1996}.}
    \label{fig:rede_reatores}
\end{figure}

The reactor network design problem involves two continuous stirred tank reactor (CSTRs), where the reaction $ \textrm{A} \rightarrow \textrm{B} \rightarrow \textrm{C} $ takes place, such that $ \textrm{A} $, $ \textrm{B} $ and $ \textrm{C} $ are the system components. In both reactions, first-order kinetics is assumed and reactors are ideally mixed, that is, there are no concentration gradients and the reactor concentration is the same as the outlet concentration. The objective is to design a system so that the concentration of $ \textrm{B} $ in the output current of the second reactor $ \left( c_{_{\mathrm{B} 2}} \right) $ is maximized, and the investment cost does not exceed a prescribed upper limit~\cite{bib:smith1996}. Figure~\ref{fig:rede_reatores} illustrates the essential structure of the reactor network design.

Let $ V_{1} $ and $ V_{2} $ be the residence times for the first and second reactors, respectively, the input current of the first reactor $ c_{_{\mathrm{A} 1}} = 1.0 \; \textrm{mol/l} $, and $ k_{1} = 0.09755988 \; \mathrm{s}^{-1} $, $ k_{2} = 0.09658428 \; \mathrm{s}^{-1} $, $ k_{3} = 0.0391908 \; \mathrm{s}^{-1} $ and $ k_{4} = 0.03527172 \; \mathrm{s}^{-1} $ the reaction rate constants, whose relationship is illustrated in Fig.~\ref{fig:rede_reatores}. In this model, it is assumed that the operating cost of a reactor is proportional to the square root of its residence time~\cite{bib:smith1996}. The reactor network design problem, whose purpose is to maximize $ f\left( c_{_{\mathrm{B} 2}} \right) = c_{_{\mathrm{B} 2}} $, is formulated as
\begin{equation} \label{eq:prob_rede_reatores_det}
\begin{split}
\max \; &f\left( c_{_{\mathrm{B} 2}} \right)\\
\textrm{Subject to} \; &c_{_{\mathrm{A} 1}} - c_{_{\mathrm{A} 0}} + k_{1} c_{_{\mathrm{A} 1}} V_{1} = 0\\
\phantom{Subject to} \; &c_{_{\mathrm{A} 2}} - c_{_{\mathrm{A} 1}} + k_{2} c_{_{\mathrm{A} 2}} V_{2} = 0\\
\phantom{Subject to} \; &c_{_{\mathrm{B} 1}} + c_{_{\mathrm{A} 1}} - c_{_{\mathrm{A} 0}} + k_{3} c_{_{\mathrm{B} 1}} V_{1} = 0\\
\phantom{Subject to} \; &c_{_{\mathrm{B} 2}} - c_{_{\mathrm{B} 1}} + c_{_{\mathrm{A} 2}} - c_{_{\mathrm{A} 1}} + k_{4} c_{_{\mathrm{B} 2}} V_{2} = 0\\
\phantom{Subject to} \; &\sqrt{V_{1}} + \sqrt{V_{2}} \leq 4\\
\phantom{Subject to} \; &0 \leq c_{_{\mathrm{A} 1}}, \; c_{_{\mathrm{A} 2}}, \; c_{_{\mathrm{B} 1}}, \; c_{_{\mathrm{B} 2}} \leq 1\\
\phantom{Subject to} \; &0 \leq V_{1}, \; V_{2} \leq 16 \; .
\end{split}
\end{equation}

This model was proposed by Manousiouthakis and Sourlas~\cite{bib:manousiouthakis1992} and, in the context of single-objective optimization, it was also studied by Ryoo and Sahinidis~\cite{bib:ryoo1995}, Smith~\cite{bib:smith1996}, Angira and Babu~\cite{bib:angira2002}, Kheawhom~\cite{bib:kheawhom2010}, Sharma~et~al.~\cite{bib:sharma2013}, and Dong~et~al.~\cite{bib:dong2014}. All authors considered the problem difficult because it has two local optima very close to the global optimum (see Table~\ref{tab:otimos_rede_reatores}). Note that both local optima use only a single reactor, while the global solution uses both. In addition, the problem contains a set of non-convex constraints, which can be challenging from the point of view of the optimization procedure.
\begin{table}[!hb]
\centering
\caption{Feasible solutions for the reactor network design problem, with an indication of the type of optimal solution (local or global) for each point~\cite{bib:ryoo1995}. Concentrations are shown in $ \left[ \textrm{mol/l} \right] $ and residence times in $ \left[ \textrm{s} \right] $.}
\label{tab:otimos_rede_reatores}
{\setlength{\tabulinesep}{1.3mm}
\setlength{\tabcolsep}{1.3mm}
\begin{tabu}{cccccccc}
\hline\hline
Type & $ c_{_{\mathrm{A} 1}} $ & $ c_{_{\mathrm{A} 2}} $ & $ c_{_{\mathrm{B} 1}} $ & $ c_{_{\mathrm{B} 2}} $ & $ V_{1} $ & $ V_{2} $ & $ f $ \\ \hline
Local & $ 0.390 $ & $ 0.390 $ & $ 0.375 $ & $ 0.375 $ & $ 16 $ & $ 0 $ & $ 0.375 $ \\
Local & $ 1 $ & $ 0.393 $ & $ 0 $ & $ 0.388 $ & $ 0 $ & $ 16 $ & $ 0.388 $ \\
Global & $ 0.771 $ & $ 0.517 $ & $ 0.204 $ & $ 0.389 $ & $ 3.037 $ & $ 5.096 $ & $ 0.389 $ \\ \hline\hline
\end{tabu}}
\end{table}

In the scenario of optimization with uncertainties, more specifically in the reliability-based optimization, Lobato~et~al.~\cite{bib:lobato2019} proposed that reasonable uncertainties would be associated with constant $ k_{1} $ to $ k_{4} $, since they are estimated using the Arrhenius temperature dependence, whose formulation involves uncertainties. The authors suggest an algebraic simplification, by eliminating equality constraints. Initially, one writes $ c_{_{\mathrm{B} 2}} $ explicitly in the last equality constraint of Eq.~(\ref{eq:prob_rede_reatores_det}). Then, the remaining equality constraints are written only as a function of $ c_{_{\mathrm{A} 1}} $ and $ c_{_{\mathrm{A} 2}} $. The resulting expressions are replaced in the corresponding terms, on the right-hand side of $ c_{_{\mathrm{B} 2}} $. Therefore, the objective function is rewritten as
\begin{equation} \label{eq:fobj_simp_rede_reatores}
f \left( \mathbf{d}, \; \mathbf{X} \right) = \dfrac{\Upsilon_{2} \left( \Upsilon_{1} \left( 1 - d_{1} \right) - \left( \Upsilon_{1} + X_{3} \left( 1 - d_{1} \right) \right) \left( d_{2} - d_{1} \right) \right)}{\left( \Upsilon_{1} + X_{3} \left( 1 - d_{1} \right) \right) \left( \Upsilon_{2} + X_{4} \left( d_{1} - d_{2} \right) \right)} \; ,
\end{equation}

\noindent with $ \Upsilon_{1} = X_{1} d_{1} $ and $ \Upsilon_{2} = X_{2} d_{2} $, where $ \mathbf{d} $ is the vector of deterministic decision variables, given by $ \mathbf{d} = \left( c_{_{\mathrm{A} 1}}, \; c_{_{\mathrm{A} 2}} \right) $, and $ \mathbf{X} = \left( X_{1}, \; X_{2}, \; X_{3}, \; X_{4} \right) $ is a set of independent random variables, with normal distribution. Thus, the problem has a single inequality constraint that, in general, is easier to manipulate, and the objective function is defined by two decision variables, $ c_{_{\mathrm{A} 1}} $ and $ c_{_{\mathrm{A} 2}} $.

In the case of the inequality constraint of Eq.~(\ref{eq:prob_rede_reatores_det}), the simplification proposed by Lobato~et~al.~\cite{bib:lobato2019} follows the same procedure previously adopted. By explicitly writing out $ V_{1} $ and $ V_{2} $ in the first two equality constraints, it becomes possible to write the inequality constraint also according to $ c_{_{\mathrm{A} 1}} $ and $ c_{_{\mathrm{A} 2}} $. Thus, an equivalent inequality constraint is given by
\begin{equation} \label{eq:rest_des_simp_rede_reatores}
\sqrt{\dfrac{1 - d_{1}}{X_{1} d_{1}}} + \sqrt{\dfrac{d_{1} - d_{2}}{X_{2} d_{2}}} - 4 \leq 0 \; .
\end{equation}

\noindent The reliability-based multi-objective optimization problem is given by
\begin{equation} \label{eq:prob_rede_reatores_conf}
\begin{split}
\max \; &\left( f \left( \mathbf{d}, \; \mathbf{X} \right), \; \beta_{\mathrm{t}} \right)\\
\textrm{Subject to} \; &\mathrm{P} \left[ g_{1} \left( \mathbf{d}, \; \mathbf{X} \right) \leq 0 \right] \leq \Phi \left( -\beta_{\mathrm{t}} \right)\\
\phantom{Subject to} \; &10^{-5} \leq d_{1}, \; d_{2} \leq 1\\
\phantom{Subject to} \; &0.1 \leq \beta_{\mathrm{t}} \leq 5 \; ,
\end{split}
\end{equation}

\noindent with $ f \left( \mathbf{d}, \; \mathbf{X} \right) $ given by Eq.~(\ref{eq:fobj_simp_rede_reatores}), and $ g_{1} \left( \mathbf{d}, \; \mathbf{X} \right) $ is the left-hand side of Eq.~(\ref{eq:rest_des_simp_rede_reatores}). For this problem, $ \mu_{i} = k_{i} $ and the standard deviation is defined in terms of the coefficient of variation, with $ \sigma_{i} / \mu_{i} = 0.15 $, for $ i = 1, \; \dots, \; 4 $. Means and standard deviations of $ \mathbf{X} $ are explicitly shown in Table~\ref{tab:prop_estat_rede_reatores}.

\begin{table}[!hb]
\centering
\caption{Mean $ \left( \vect{\mu} \right) $ and standard deviation $ \left( \vect{\sigma} \right) $ of the random variables of the reactor network design problem~\cite{bib:lobato2019}.}
\label{tab:prop_estat_rede_reatores}
{\setlength{\tabulinesep}{1.3mm}
\setlength{\tabcolsep}{1.3mm}
\begin{tabu}{lcccc}
\hline\hline
& $ X_{1} $ & $ X_{2} $ & $ X_{3} $ & $ X_{4} $ \\ \hline
$ \vect{\mu} $ & $ 0.09755988 $ & $ 0.09658428 $ & $ 0.0391908 $ & $ 0.03527172 $ \\
$ \vect{\sigma} $ & $ 0.01463398 $ & $ 0.01448764 $ & $ 0.00587862 $ & $ 0.00529076 $ \\ \hline\hline
\end{tabu}}
\end{table}

The problem is solved using two different approaches. At first, Eq.~(\ref{eq:prob_rede_reatores_conf}) is optimized using MODE, with the same control parameters as in previous executions (Sections~\ref{sec:benchmark_novel_approach} and \ref{sec:heat_exchanger_novel_approach}). In the second approach, the problem of Eq.~(\ref{eq:prob_rede_reatores_conf}) is converted into a single-objective problem, whose only objective is given by $ f \left( \mathbf{d}, \; \mathbf{X} \right) $. Thus, 50 independent executions of DE are performed, varying the reliability index, which is defined by equally spaced values in the interval $ 0.1 \leq \beta_{\mathrm{t}} \leq 5 $. Each search process terminates when the number of generations exceeds 100. In both approaches, the probabilistic constraints are solved using ASOSL with control parameters as in Sections~\ref{sec:benchmark_novel_approach} and \ref{sec:heat_exchanger_novel_approach}.

\begin{figure}[!ht]
    \centering
    \includegraphics[width=\columnwidth]{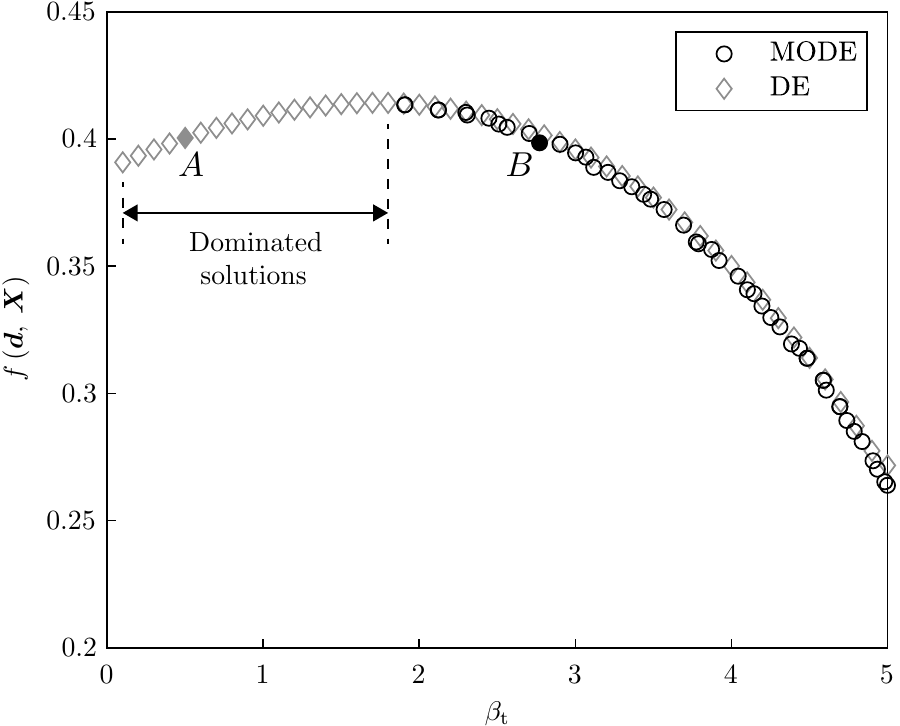}
    \caption{Optimal results obtained by MODE for a reactor network design formulated as a multi-objective optimization problem given by Eq.~(\ref{eq:prob_rede_reatores_conf}), and by DE for an equivalent single-objective problem, varying the reliability index in each execution. Here, the objectives $ f $ and $ \beta_{\mathrm{t}} $ are to be maximized. Points $ A $ and $ B $, represented by the filled diamond and circle, respectively, have approximately the same value for $ f $. However, the reliability index in $ B $ is higher. Therefore, $ B $ dominates $ A $.}
    \label{fig:MODE_DE_reactor_network}
\end{figure}

Figure~\ref{fig:MODE_DE_reactor_network} presents the results obtained by both methods on the same axes. From this perspective, the relevance of modeling the reliability-based optimization problem in a multi-objective context is clear: in the interval analyzed for the reliability index, DE obtains several dominated solutions (which are useless), since dominance criteria---in Pareto sense---are not taken into account (this holds for most metaheuristics in the single-objective context). DE is not able to identify whether the reliability index is adequate, as long as each execution is independent, but only to maximize the objective function, for a given $ \beta_{\mathrm{t}} $ fixed. All solutions obtained by DE when $ \beta_{\mathrm{t}} < 1.8 $ are dominated by another point, with equal value for $ f $ but higher reliability index---they do not belong to the Pareto set.

In turn, the multi-objective formulation of Eq.~(\ref{eq:prob_rede_reatores_conf}) allows the metaheuristic to obtain a set of results that is composed only of dominant solutions. Commonly, this fact is also reflected in the computational cost. The single-objective problem must be solved as many times as the number of solutions desired (according to the number of different reliability indices), while the multi-objective approach provides a set of dominant solutions with a single execution. As a result, the multi-objective approach does not compute useless solutions.

Now, consider the problem given by Eq.~(\ref{eq:prob_rede_reatores_conf}) and the incorporation of robustness analysis. The formulation of the reliability-based robust multi-objective optimization problem is defined as
\begin{equation} \label{eq:prob_rede_reatores_conf_rob}
\begin{split}
\max \; &\left( f \left( \mathbf{d}, \; \mathbf{X}, \; \vect{\delta} \right), \; \beta_{\mathrm{t}} \right)\\
\textrm{Subject to} \; &\mathrm{P} \left[ g_{1} \left( \mathbf{d}, \; \mathbf{X}, \; \vect{\delta} \right) \leq 0 \right] \leq \Phi \left( -\beta_{\mathrm{t}} \right)\\
\phantom{Subject to} \; &10^{-5} \leq d_{1}, \; d_{2} \leq 1\\
\phantom{Subject to} \; &0.1 \leq \beta_{\mathrm{t}} \leq 5 \; .
\end{split}
\end{equation}

\begin{figure}[!ht]
    \centering
    \includegraphics[width=\columnwidth]{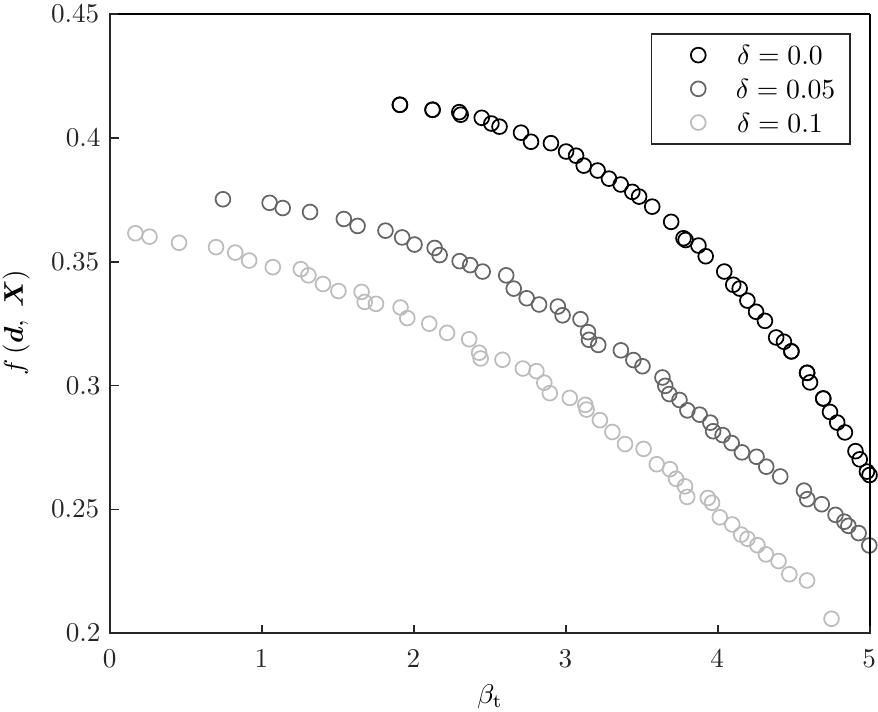}
    \caption{Results of the multi-objective optimization problem with uncertainty of Eq.~(\ref{eq:prob_rede_reatores_conf_rob}), for different levels of robustness.}
    \label{fig:reactor_network_rob_conf}
\end{figure}

\noindent In this case, it is proposed to evaluate the robustness of $ c_{_{\mathrm{A} 1}} $ and $ c_{_{\mathrm{A} 2}} $. For this purpose, the effective mean technique is used, with $ M = 50 $ and noise set at 5\% and 10\%. Figure~\ref{fig:reactor_network_rob_conf} shows the profiles obtained for both values of $ \delta $, in addition to the deterministic case $ \left( \delta = 0 \right) $. Again, the results are obtained using MODE and ASOSL, with the same parameters as in the previous cases.

As in the previous cases, there is a shift in the profiles as the value of $ \delta $ is increased. However, this model has three particularities. First, the Pareto front curvature (the simulations depicted in Fig.~\ref{fig:reactor_network_rob_conf} are a discretization of the Pareto curve) vary with the noise value $ \delta $, becoming flatter when $ \delta $ increases. Second, for a given reliability index, this is reflected in a non-linear reduction in the objective values. As an example, let $ \beta_{\mathrm{t}} = 3 $. The decrease in the value of $ f $ is more pronounced between the values calculated when $ \delta = 0 $ and $ \delta = 0.05 $, than for $ \delta = 0.05 $ and $ \delta = 0.1 $. This fact must be taken into account when choosing the result to be implemented in practice. The third particularity concerns the influence of robustness on the reliability of the results. For a given target value of objective $ f $, increasing $ \delta $ means giving up of more reliable optima, in terms of the reliability index. Consider, for example, $ f = 0.3 $ in Figure~\ref{fig:reactor_network_rob_conf}. The displacement of Pareto curves with different values for $ \delta $ leads to the need for a decision between choices that favor either the robustness or the reliability of the results. Thus robustness and reliability, measured by $ \delta $ and $ \beta_{\mathrm{t}} $, respectively, are also conflicting objectives.

Another relevant issue concerns the dispersion of each Pareto curve. This outcome, observed in Figs.~\ref{fig:rede_trocadores_pareto} and \ref{fig:reactor_network_rob_conf}, can be further explored by adjusting the result, obtained for each $ \delta $, by a polynomial of degree two. Thus, it is possible to quantify the deviations in relation to the adjusted curves, to estimate the error made, and the dispersion of the set of points. Table~\ref{tab:erro_ajuste_resultados} shows some goodness-of-fit statistics for adjusted data: (\textit{i}) the sum of squares due to error (SQR); (\textit{ii}) R-square $ \left( R^{2} \right) $; (\textit{iii}) adjusted R-square $ \left( R_{\mathrm{a}}^{2} \right) $ and; (\textit{iv}) root mean squared error (RMS)~\cite{bib:gujarati2008}.

\begin{table}[!ht]
\centering
\caption{Goodness-of-fit statistics for each level of robustness in Fig.~\ref{fig:reactor_network_rob_conf}.}
\label{tab:erro_ajuste_resultados}
{\setlength{\tabulinesep}{1.3mm}
\setlength{\tabcolsep}{1.5mm}
\begin{tabu}{lcccc}
\hline\hline
 & SQR & $ R^{2} $ & $ R_{\mathrm{a}}^{2} $ & RMS \\ \hline
$ \delta = 0 $ & $ 9.483 \times 10^{-5} $ & $ 0.9991 $ & $ 0.9991 $ & $ 1.452 \times 10^{-3} $ \\
$ \delta = 0.05 $ & $ 2.851 \times 10^{-4} $ & $ 0.9969 $ & $ 0.9968 $ & $ 2.463 \times 10^{-3} $ \\
$ \delta = 0.1 $ & $ 3.828 \times 10^{-3} $ & $ 0.9604 $ & $ 0.9596 $ & $ 8.93 \times 10^{-3} $ \\ \hline\hline
\end{tabu}}
\end{table}

Analyzing the statistics, the dispersion of the set of points tends to increase when the value of $ \delta $ increases. This behavior may be due to the randomness inserted by the robustness analysis, using the effective mean technique. The increase of $ \delta $ is related to the interval in the neighborhood of a candidate solution, where random points are spread. Therefore, the average value of the objective, calculated by using these random points, influences the dispersion of the approximation of the Pareto curves by the metaheuristic.

\begin{figure}[!ht]
    \centering
    \includegraphics[width=\columnwidth]{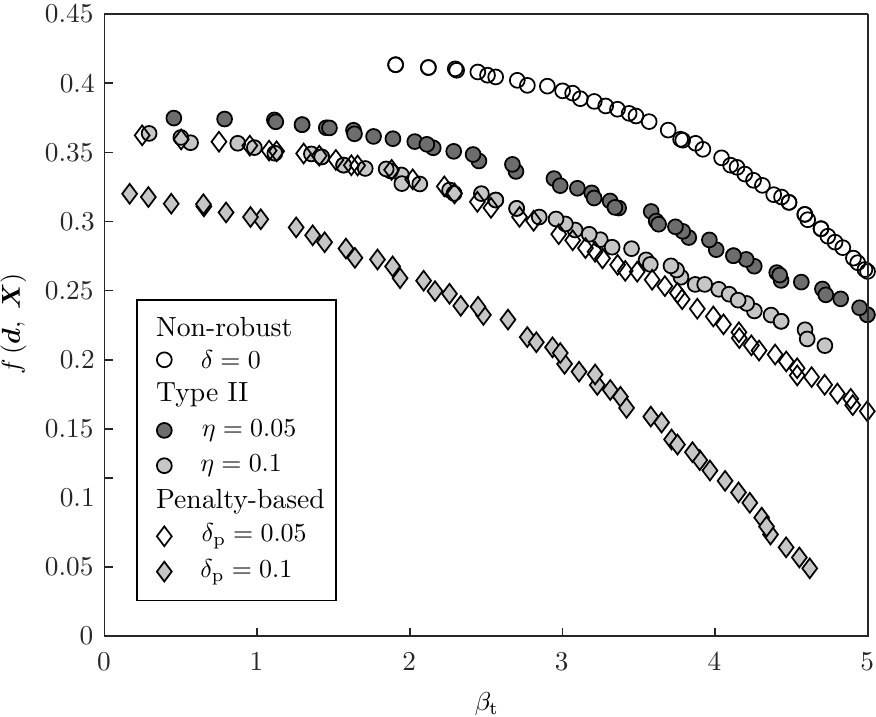}
    \caption{Pareto curves obtained using Type II and Penalty-based approaches, for robustness levels set at 0\%, 5\% and 10\%.}
    \label{fig:reactor_network_rob_conf_penal_moddiff}
\end{figure}

In view of this scenario, it is pertinent to obtain the results of the problem given by Eq.~(\ref{eq:prob_rede_reatores_conf_rob}), using other techniques for robustness analysis. Figure~\ref{fig:reactor_network_rob_conf_penal_moddiff} shows the Pareto curves obtained using Type II and Penalty-based approaches (discussed in Sections~\ref{sec:tec_robustez} and \ref{sec:tec_rob_penal}, respectively), for robustness levels set at 5\% and 10\%, in addition to the deterministic case. In all runs, the results were computed using MODE and ASOSL, with the same values for control parameters adopted in previous executions.

Initially, note that none of the approaches is capable of reaching a reliability index close to $ \beta_{\mathrm{t}} = 5 $ when the robustness is specified at $ 10 \% $, as in the case of effective mean strategy (see Fig.~\ref{fig:reactor_network_rob_conf}). In addition, it is clear that the Penalty-based technique introduces a more severe penalty in relation to other techniques. In turn, Type II technique presents results very close to those achieved using effective mean. Regarding the dispersion of the results, it is clear that this is not a particular issue of the effective mean technique, since it turns up in all profiles shown in Fig.~\ref{fig:reactor_network_rob_conf_penal_moddiff}.

\subsection{Catalyst mixing problem}

\subsubsection{Single-objective formulation}

In this problem, an optimal policy for mixing two catalysts along the length of a tubular reactor must be determined. The problem involves the reactions
\begin{equation*}
S_{1} \xleftrightarrow{\hspace{0.25em} k_{1}, \; k_{2} \hspace{0.25em}} S_{2} \xrightarrow[]{\hspace{0.25em} k_{3} \hspace{0.25 em}} S_{3} \; ,
\end{equation*}

\noindent where $ k_{1} $ and $ k_{2} $ are the reaction rate constants of the first two reactions, in a reactor in which the catalyst consists essentially of the compound that catalyzes the reversible reactions $ S_{1}~\longleftrightarrow~S_{2} $, and $ k_{3} $ is the reaction rate constant of $ S_{2}~\longrightarrow~S_{3}~$. The model, formulated by Gunn and Thomas~\cite{bib:gunn1965}, aims to maximize the production of $ S_{3} $. Both catalysts activate the three reactions, but the activation rates are different for each reaction pair and catalyst agents. Therefore, different catalyst mixing policies leads to a wide range of product compositions.

Let
\begin{subequations} \label{eq:prob_geral_mistura_catalisador}
\begin{equation} \label{eq:objetivo_problema_catalisadores}
f \left[ v \right] = 1 - y_{1} \left( t_{\mathrm{f}} \right) - y_{2} \left( t_{\mathrm{f}} \right) \; ,
\end{equation}
\noindent be the objective function for the catalyst mixing problem defined as
\begin{equation} \label{eq:prob_mistura_catalisador}
\begin{split}
\max \; &f \left[ v \right]\\
\textrm{Subject to} \; &\dfrac{\mathrm{d} y_{1}}{\mathrm{d} t} = v \left( t \right) \left( 10 y_{2} \left( t \right) - y_{1} \left( t \right) \right)\\
\phantom{Subject to} \; &\dfrac{\mathrm{d} y_{2}}{\mathrm{d} t} = v \left( t \right) \left( y_{1} \left( t \right) - 10 y_{2} \left( t \right) \right) - \left( 1 - v \left( t \right) \right) y_{2} \left( t \right)\\
\phantom{Subject to} \; &y_{1} \left( 0 \right) = 1, \; y_{2} \left( 0 \right) = 0\\
\phantom{Subject to} \; &0 \leq v \left( t \right) \leq 1\\
\phantom{Subject to} \; &0 \leq t \leq t_{\mathrm{f}}\; ,
\end{split}
\end{equation}
\end{subequations}

\noindent where $ y_{1} $ and $ y_{2} $ are, respectively, the molar fractions of the substances $ S_{1} $ and $ S_{2} $, and $ t $ represents the residence time of the substances in the reactor, from the time of entry into the reactor $ t_{0} $ (here taken as zero), to the time of exit $ t_{\mathrm{f}} $. The design variable, in this problem, is a control variable $ v \left( t \right) $, which represents the catalyst mixing fraction, a portion of the catalyst obtained from the substance that catalyzes the reaction $ S_{1}~\longleftrightarrow~S_{2} $, which must be determined by the optimizer~\cite{bib:vassiliadis1993,bib:dadebo1995,bib:tanartkit1997,bib:bell2000,bib:irizarry2005,bib:liu2013,bib:liu2015,bib:lobato2019}.

The problem is solved by discretizing the control variable $ v \left( t \right) $ and the time variable. For this purpose, time is discretized and, for each subinterval $ \left[ t_{k}, \; t_{k + 1} \right] $, a fraction of the catalyst is approximated and is represented by $ v_{k} $, for $ k = 0, \; \dots, \; n_{t} - 1 $, where $ n_{t} $ represents the number of time subintervals. Therefore, $ \mathbf{d} = \left( v_{0}, \; \dots, \; v_{n_{t} \, - \, 1}, \; t_{1}, \; \dots, \; t_{n_{t} \, - \, 1} \right) $ is the decision vector that must be calculated by the optimizer.

Specifically for the case where $ n_{t} = 3 $, the optimization algorithm must determine $ v_{k} $, for $ k = 0, \; 1, \; 2 $, in addition to the intermediate time instants, $ t_{1} $ and $ t_{2} $ ($ t_{0} = 0 $ and $ t_{3} = t_{\mathrm{f}} = 1 $ are known). Using DE with the same values adopted previously for the control parameters, the following results are obtained, after 100 generations: $ v_{0} = 1 $ for $ t \leq 0.1338 $, $ v_{1} = 0.2248 $ for $ 0.1338 < t \leq 0.7237 $, and $ v_{2} = 0 $ for $ 0.7237 < t \leq t_{\mathrm{f}} $. Under these conditions, the optimum value of the objective function is $ f = 0.048065 $. This result agrees with Gunn and Thomas~\cite{bib:gunn1965}, Vassiliadis~\cite{bib:vassiliadis1993}, and Lobato~et~al.~\cite{bib:lobato2011b}.

\subsubsection{Reliability and robustness analysis}

Lobato~et~al.~\cite{bib:lobato2019} suggest that the multi-objective optimization problem can be formulated considering potential uncertainties in terms of the discretized control variable. In this case, let $ \mathbf{X} = \left( X_{0}, \; X_{1}, \; X_{2} \right) $ be a set of independent random variables (representing the discretized control variables), where $ X_{i} \sim \mathcal{N} \left( \mu_{i}, \; \sigma_{i}^{2} \right) $, for $ i = 0, \; 1, \; 2 $. The standard deviation of $ \mathbf{X} $ is defined in relation to the coefficient of variation, $ \sigma_{i} = 0.1 \mu_{i} $, where $ \mu_{i} $ represents the mean of the random variables. In the case of robustness analysis, $ \vect{\delta} $ represents the noise vector, related to the discretized control variable $ v_{k} $, for $ k = 0, \; 1, \; 2 $. In this case, the deterministic decision vector is composed by $ \mathbf{d} = \left( t_{1}, \; t_{2} \right) $. The problem defined by Eq.~(\ref{eq:prob_geral_mistura_catalisador}), but formulated as a multi-objective optimization problem with uncertainties, is given by
\begin{equation} \label{eq:prob_mistura catalisador_incertezas}
\begin{split}
\max \; &\left( f \left( \mathbf{d}, \; \vect{\mu}, \; \vect{\delta} \right), \; \beta_{\mathrm{t}} \right)\\
\textrm{Subject to} \; &\dfrac{\mathrm{d} y_{1}}{\mathrm{d} t} = X \left( 10 y_{2} \left( t \right) - y_{1} \left( t \right) \right)\\
\phantom{Subject to} \; &\dfrac{\mathrm{d} y_{2}}{\mathrm{d} t} = X \left( y_{1} \left( t \right) - 10 y_{2} \left( t \right) \right) - \left( 1 - X \right) y_{2} \left( t \right)\\
\phantom{Subject to} \; &y_{1} \left( 0 \right) = 1, \; y_{2} \left( 0 \right) = 0\\
\phantom{Subject to} \; &0 \leq d_{1} \leq 0.5\\
\phantom{Subject to} \; &0.5 < d_{2} \leq 1\\
\phantom{Subject to} \; &\mathrm{P} \left[ g_{i} \left( \mathbf{X}, \; \vect{\delta} \right) \leq 0 \right] \leq \Phi \left( -\beta_{\mathrm{t}} \right)\\
\phantom{Subject to} \; &0 < \mu_{i} \leq 1\\
\phantom{Subject to} \; &0.1 \leq \beta_{\mathrm{t}} \leq 5 \; ,
\end{split}
\end{equation}

\noindent where $ g_{i} \left( \mathbf{X}, \; \vect{\delta} \right) = X_{i} $, for $ i = 1, \; 2, \; 3 $.

To determine the set of feasible solutions, initially one generates an initial population of candidate solutions. For each evaluation of the objective function, the vector of decision variables is perturbed, in order to evaluate robustness. Perturbations occur only in the coordinates associated to the discretized control variable (time and reliability index are not assumed susceptible to noise).

In the case of probabilistic constraints, an inverse reliability analysis method provides the optimal values of the discretized control variable, as well as the respective values of the performance function. These results, as well as the discretized time instants, are employed to obtain the solution of coupled differential equations. Thus, the molar fractions $ y_{1} \left( t \right) $ and $ y_{2} \left( t \right) $ are used to evaluate the objective given by Eq.~(\ref{eq:objetivo_problema_catalisadores}).

\begin{figure}[!ht]
    \centering
    \includegraphics[width=\columnwidth]{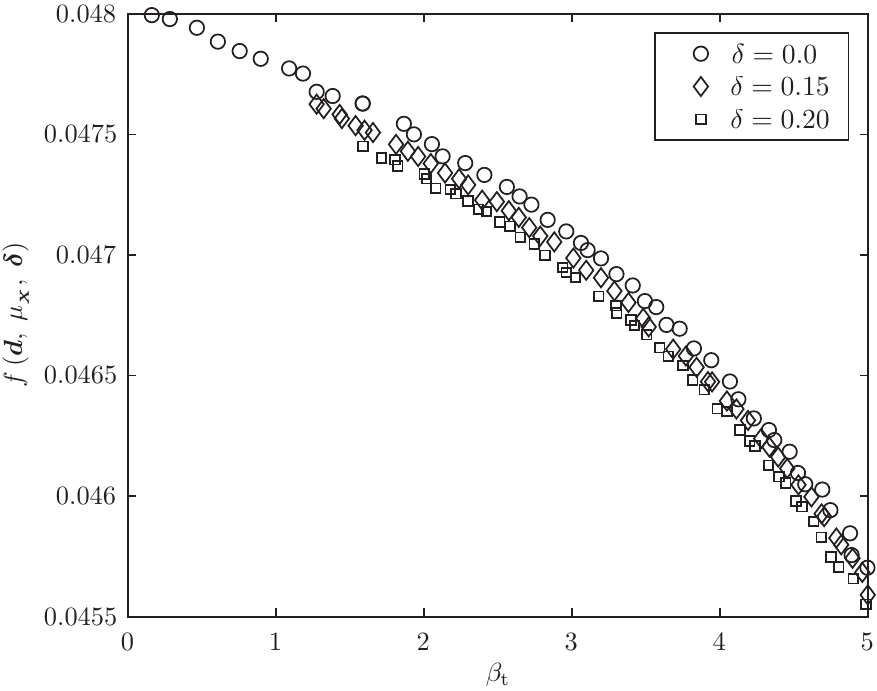}
    \caption{Results of the multi-objective catalyst mixing problem with uncertainties, given by Eq.~(\ref{eq:prob_mistura catalisador_incertezas}), for different levels of robustness.}
    \label{fig:catalyst_rob_conf}
\end{figure}

MODE algorithm is used to obtain the Pareto set to Eq.~(\ref{eq:prob_mistura catalisador_incertezas}), the probabilistic constraints are evaluated with ASOSL and the robustness analysis is performed by the effective mean strategy, all using the same parameter values adopted in the previous problems. Figure~\ref{fig:catalyst_rob_conf} presents the Pareto curves, calculated for $ \delta = 0.15 $ and $ \delta = 0.20 $, as well as for the deterministic case $ \left( \delta = 0 \right) $.

It is possible to compare the behavior of the Pareto sets obtained, in relation to the set of results in which external noise is not considered. In general, engineering problems are analyzed taking into account lower noise levels. When this optimization problem is evaluated for values of $ \delta $ closer to zero, the associated Pareto curves practically overlap with the results obtained for $ \delta = 0 $. However, it is relevant to assess this problem because it makes clear that not always the robustness is relevant, since some functions may not be very sensitive to noise.

\begin{figure}[!ht]
    \centering
    \includegraphics[width=\columnwidth]{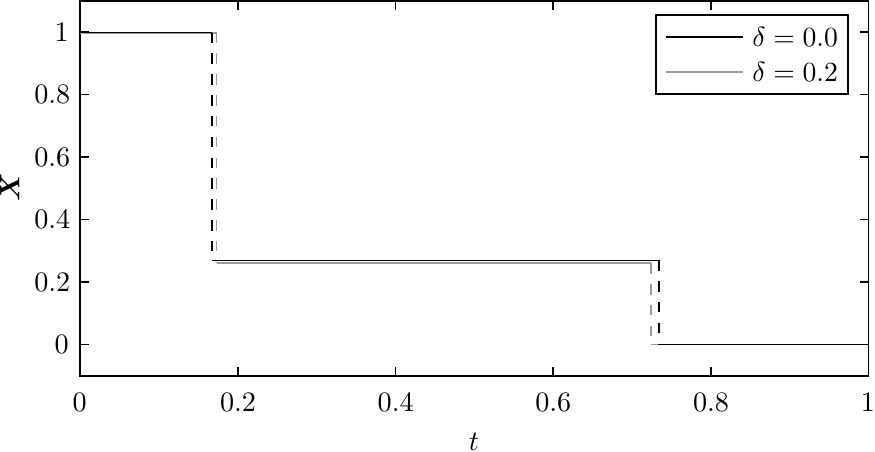}
    \caption{Control profiles obtained for two different levels of robustness, when $ \beta_{\mathrm{t}} \approx 1.6 $. There is a slight difference between the profiles, even with excessive variation in the noise level. The values of the control variables are shown in Table~\ref{tab:variaveis_controle}.}
    \label{fig:controle_delta}
\end{figure}

This behavior can be best illustrated by analyzing Fig.~\ref{fig:controle_delta}, where the control profiles of the problem given by Eq.~(\ref{eq:prob_mistura catalisador_incertezas}) are shown, when $ \delta = 0 $ and $ \delta = 0.2 $, with reliability indexes approximately equal to $ \beta_{\mathrm{t}} \approx 1.6 $ for both cases. The values of the discretized control variables are shown in Table~\ref{tab:variaveis_controle}, for each level of robustness analyzed. For such decision variables, the corresponding objectives are $ f = 0.0476 $ and $ f = 0.0474 $, respectively.

\begin{table}[!ht]
\centering
\caption{Values of the discretized control variables for two different levels of robustness $ \left( \delta \right) $, when $ \beta_{\mathrm{t}} \approx 1.6 $.}
\label{tab:variaveis_controle}
{\setlength{\tabulinesep}{1.5mm}
\setlength{\tabcolsep}{2mm}
\begin{tabu}{cccc}
\hline\hline
$ \delta $ & \multicolumn{3}{c}{Control variable} \\ \hline
\multirow{2}{*}{0.0} & $ x_{0} = 0.9984 $ & $ x_{1} = 0.2695 $ & $ x_{2} = 0.0004 $ \\
 & $ t \leq 0.1669 $ & $ 0.1669 < t < 0.7341 $ & $ t \geq 0.7341 $ \\
\multirow{2}{*}{0.2} & $ x_{0} = 0.9982 $ & $ x_{1} = 0.2601 $ & $ x_{2} = 0.0002 $ \\
 & $ t \leq 0.1728 $ & $ 0.1728 < t < 0.7242 $ & $ t \geq 0.7242 $ \\ \hline\hline
\end{tabu}}
\end{table}

In both cases, $ x_{0} \approx 1 $ and $ x_{2} \approx 0 $. Thus, the objective value is essentially determined by $ x_{1} $, and by the times of occurrence, $ t_{1} $ and $ t_{2} $. Therefore, a significant deviation from the control parameter $ x_{1} $, in relation to cases in which $ \delta = 0 $ and $ \delta = 0.2 $, could represent a substantial difference between the corresponding objective values. However, note that the control profiles obtained in the cases analyzed in Fig.~\ref{fig:controle_delta} are very similar. These profiles are even similar to the deterministic solution, shown previously. This means that the problem, formulated by Eq.~(\ref{eq:prob_mistura catalisador_incertezas}), has solutions not very sensitive to external perturbations, which corroborates the slight difference of the Pareto curves obtained in Fig.~\ref{fig:catalyst_rob_conf}, even with excessive variation in the noise level.

\section{Conclusions} \label{sec:conclusions}

The novel formulation proposed in this work, regarding a reliability-based robust design multi-objective optimization problem, presents a systematic and consistent way of considering uncertainties during the optimization procedure of computational models. The contribution of each of these approaches---robustness and reliability---to the results obtained is clear, which makes it possible to understand the properties and purpose of them, and modeling the problem according to the needs of the project. The formulation also allows considering uncertainties of different types in any of the decision variables and inequality constraints.

Although robustness and reliability are conflicting with each other, in terms of the reliability index and the noise level, the analysis of the two-dimensional benchmark problem provides a clearer understanding of what is going on. The insertion of a third objective in the formulation of the problem (possibly treating the noise level as an additional objective) would mean obtaining a Pareto front, in the objective space, represented by a surface, which may be difficult to analyze.

A certain dispersion of the optimal values was observed in all problems analyzed, when the noise levels of the decision variables are high. This is not exactly related to the robustness analysis method, nor to the problem formulation, but to the strategy used to evaluate the sensitivity of a given candidate solution, taking into account a set of random points. The metaheuristic evaluates the dominance of each individual in the population, in order to define the fittest ones, considering the average contribution of the random samples. Consequently, the sensitivity of the objective function at a given point affects the genetic operators of the optimization algorithm, causing the profiles obtained to show such dispersion.

In general, the results are able to demonstrate the importance of considering uncertainties during the optimization procedure. It is noted that robust and reliable solutions tend to deviate from deterministic optima, which corroborates the fact that deterministic optimizers can hardly be implemented in practice, since they may be very sensitive to perturbations, which can prevent them from meeting the needs of the system when confronted with real operating conditions.



\section*{Acknowledgments}

\noindent This study was financed in part by the Coordination for the Improvement of Higher Education Personnel---Brasil (CAPES)---Finance Code 001. Gustavo Libotte is supported by a postdoctoral fellowship from the Institutional Training Program (PCI) of the Brazilian National Council for Scientific and Technological Development (CNPq), grant number 303185/2020-1, and Carlos Chagas Filho Foundation for Supporting Research in the State of Rio de Janeiro, grant number E-26/200.560/2018.

\bibliographystyle{elsarticle-num} 
\bibliography{bibliography}
\end{document}